  \documentclass[12pt]{amsart}
  \usepackage{latexsym} 
  \usepackage[all]{xy}
  \usepackage{amsfonts} 
  \usepackage{amsthm} 
  \usepackage{amsmath} 
  \usepackage{amssymb}
  \usepackage{pifont} 
  \newcommand{\cC}{{C}}  
\numberwithin{equation}{section}

  \def\sw#1{{\sb{(#1)}}} 
  \def\sco#1{{\sb{[#1]}}} 
  \def\su#1{{\sp{[#1]}}}  
  \def\suc#1{{\sp{(#1)}}} 
  \def\sut#1{{\sp{\langle #1\rangle}}}

  \def\<{{\langle}} 
  \def\>{{\rangle}}

  \def\eps{\varepsilon}

  \def\note#1{{}} 
  \def\can{{\rm can}}

  \def\note#1{} 
  \def\M{{\bf M}}

  \def\cM{{\mathcal M}}

  \def\cC{{\mathcal C}} 
   
  \def\cD{{\mathcal D}}

  \def\cQ{{\mathcal Q}}

  \def\can{{\rm can}}

  \def\beq{\begin{equation}} 
  \def\eeq{\end{equation}} 
  \def\DC{{\Delta_\cC}} 
  \def \eC{{\eps_\cC}}

  \def\im{{\rm Im}}

  \def\ot{{\otimes}}

  \def\roM{\varrho^{M}}

  \newcommand{\Ra}{\Rightarrow}

  \def\roA{{\varrho^A}}
   \def\Tro{{{}^T\!\varrho}}
    \def\Mro{{{}^M\!\varrho}}
  \def\roT{{\varrho^T}}
  \def\roM{{\varrho^M}}

\def\hH{\mathcal H}

  \newcounter{zlist} 
  \newenvironment{zlist}{\begin{list}{(\arabic{zlist})}{ 
  \usecounter{zlist}\leftmargin2.5em\labelwidth2em\labelsep0.5em 
  \topsep0.6ex
  \parsep0.3ex plus0.2ex minus0.1ex}}{\end{list}}

  \newcounter{blist} 
  \newenvironment{blist}{\begin{list}{(\alph{blist})}{ 
  \usecounter{blist}\leftmargin2.5em\labelwidth2em\labelsep0.5em 
  \topsep0.6ex 
  \parsep0.3ex plus0.2ex minus0.1ex}}{\end{list}} 

  \newcounter{rlist} 
  \newenvironment{rlist}{\begin{list}{(\roman{rlist})}{ 
  \usecounter{rlist}\leftmargin2.5em\labelwidth2em\labelsep0.5em 
  \topsep0.6ex 
  \parsep0.3ex plus0.2ex minus0.1ex}}{\end{list}}

\def\stac#1{\raise-.2cm\hbox{$\,\stackrel{\displaystyle\otimes}{\scriptscriptstyle{#1}}\,$}}
\def\sstac#1{\otimes_{#1}}
\def\cten#1{\raise-.2cm\hbox{$\,\stackrel{\displaystyle\widehat{\otimes}}
{\scriptscriptstyle{#1}}\,$}}

  \def\Label#1{\label{#1}\ifmmode\llap{[#1] }\else 
  \marginpar{\smash{\hbox{\tiny [#1]}}}\fi} 
  \def\Label{\label}

  \newtheorem{proposition}{Proposition}[section] 
  \newtheorem{lemma}[proposition]{Lemma} 
   
  \newtheorem{corollary}[proposition]{Corollary} 
  \newtheorem{theorem}[proposition]{Theorem} 

  \theoremstyle{definition} 
  \newtheorem{definition}[proposition]{Definition} 
  \newtheorem{example}[proposition]{Example} 

  \theoremstyle{remark} 
  \newtheorem{remark}[proposition]{Remark}

  \newcounter{c} 
   
  \newcommand{\etyk}[1]{\vspace{-7.4mm}$$\begin{equation}\Label{#1} 
  \addtocounter{c}{1}} 
  \renewcommand{\]}{\ifnum \value{c}=1 $$\else \end{equation}\fi} 
\def\cstac#1{\raise-.2cm\hbox{$\,\stackrel{\displaystyle\Box}{\scriptscriptstyle{#1}}\,$}}

\begin{document} 

 \title{Pre-torsors and equivalences} 
 \author{Gabriella B\"ohm}
 \address{Research Institute for Particle and Nuclear Physics, Budapest, 
 \newline\indent H-1525
 Budapest 114, P.O.B.\ 49, Hungary}
  \email{G.Bohm@rmki.kfki.hu}
  \author{Tomasz Brzezi\'nski}    
  \address{ Department of Mathematics, Swansea University, 
  Singleton Park, \newline\indent  Swansea SA2 8PP, U.K.} 
 \email{T.Brzezinski@swansea.ac.uk}   
    \date{July 2006; revised April 2007} 
  \begin{abstract} 
 Properties of (most general) non-commutative torsors or $A$-$B$ torsors are 
 analysed. Starting with pre-torsors it is shown that they are equivalent to a
 certain  class of Galois extensions of algebras by corings. It is shown that
 a class of faithfully flat pre-torsors induces equivalences between
 categories of comodules of associated corings. It is then proven that $A$-$B$
 torsors  correspond to monoidal functors (and, under some additional
 conditions, equivalences) between categories of  comodules of  bialgebroids. 
    \end{abstract} 
  \maketitle 
  \section{Introduction}
The notion of an $A$-$B$ torsor appeared in algebra as a result of a chain of 
natural generalisations of the notion of a quantum torsor introduced in
\cite{Gru:tor} as a formalisation of the proposal made by Kontsevich 
in \cite[Section~4.2]{Kon:ope}. As observed in \cite{Sch:tor}, a faithfully
flat quantum torsor is the same as a faithfully flat Hopf-(bi)Galois
object. To cover the case of a Hopf-Galois extension (rather than just an
object), the notion of a $B$ torsor was introduced in  \cite{Sch:tor} (cf.\
\cite[Section~2.8]{Sch:Hop}). Hence $B$ torsors can be understood as
Hopf-Galois extensions {\em without} the explicit mention of a Hopf
algebra. Furthermore, a faithfully flat $B$ torsor corresponds not only to a
Hopf-Galois extension of the base algebra $B$ on one side, but also to a
Galois object (by a $B$-bialgebroid) on the other. In order to remove this
asymmetry, the notion of an $A$-$B$ torsor was introduced in
\cite[Chapter~5]{Hobst:PhD}. A  faithfully flat $A$-$B$ torsor can be
understood as a bi-Galois (i.e.\ two-sided Galois) extension by bialgebroids.

As observed in \cite[Theorem~5.2.10]{Hobst:PhD},
 with every faithfully flat $A$-$B$ torsor $T$ there are
associated two bialgebroids, one over $A$, the other over $B$. 
These bialgebroids coact (freely) on the torsor,
making it a bicomodule. Thus $A$-$B$ torsors are natural objects which can 
facilitate a description of (monoidal) functors between categories of
 comodules of bialgebroids. 
The freeness of coactions of bialgebroids on $T$ (understood as the {\em
Galois condition}) has a very natural geometric interpretation. 
Recall that a principal bundle over a Lie groupoid is a manifold with a free 
groupoid action, whose fixed manifold coincides with the base manifold of
another Lie groupoid (cf.\ \cite[Section~5.7]{MoeMrc:fol}). Bialgebroids can
 be seen as  {\em groupoids} in non-commutative geometry (and 
for this reason are often referred to as {\em quantum groupoids}). Thus
faithfully flat $A$-$B$ torsors can be understood as {\em quantum  principal
bundles over quantum groupoids}. 
Since a monoidal functor between the comodule categories of two bialgebroids
maps comodule algebras to comodule algebras, such functors play an important
role when Galois extensions by different bialgebroids need to be related. 
This is the case, for example, in the reduction of a quantum principal 
bundle over a quantum groupoid $\cC$, to a bundle whose structure is
described by a quotient of $\cC$. 

The aim of this paper is to analyse algebraic properties of {$A$-$B$
  torsors} and their generalisation, termed pre-torsors. We start in 
Section~\ref{sec.prelims} with
recalling some preliminary results about corings, bialgebroids and
$\times_A$-Hopf algebras. The notion of an
$A$-$B$ pre-torsor is introduced in Section~\ref{sec.pre-tors}. It is shown
that faithfully 
flat pre-torsors are in bijective correspondence with faithfully flat
coring-Galois 
extensions. In particular every faithfully flat $A$-$B$ pre-torsor induces an
$A$-coring $\cC$ and a $B$-coring $\cD$. 
In the first part of 
Section~\ref{sec.equpretor} it is shown that both corings $\cC$ and $\cD$ 
arising from a faithfully flat pre-torsor $T$ can be
identified with cotensor products (as bicomodules). If $T$ satisfies some
additional 
conditions (cf.\ Remark~\ref{rem:ff} (iii)), then it induces an equivalence
between the 
categories of (left) comodules of $\cC$ and $\cD$. It is furthermore shown that
these additional conditions reduce to a natural faithful flatness assumption
provided the entwining maps induced by $T$ are
bijective. Section~\ref{sec.tors} 
deals with properties specific to $A$-$B$ torsors. Following \cite{Hobst:PhD}
we establish a bijective correspondence between faithfully flat $A$-$B$
torsors and faithfully flat Galois extensions with $\times_A$- or
$\times_B$-Hopf algebras. 
In the case when, for an algebra $B$ and an $A$-bialgebroid $\cC$, the functor
from the category of $\cC$-comodules to the category of $B$-$B$ bimodules, 
induced by a $B^e$-$\cC$ bicomodule $T$, preserves colimits, we establish a
bijective correspondence between $B$-ring and $\cC$-comodule algebra structures
in $T$ on one hand and lax monoidal structures of the induced functor on the
other hand.  
The question when is this monoidal structure strict is addressed. In
particular, it is proven that a faithfully flat $A$-$B$ torsor which is also 
faithfully flat as a right $B$-module induces a strict monoidal functor. 
Finally a class of $A$-$B$ torsors inducing monoidal equivalences 
between categories of comodules
is found. The paper is concluded with  two appendices.  
In the first one a $\times_B$-Hopf
algebra corresponding to a torsor coming from a cleft extension by a Hopf
algebroid 
is computed. This turns out to be a generalisation of bialgebroids studied in 
 \cite{ConMos:Ran} and \cite{Kad:Pseudo} (and shown to be mutually isomorphic
 in   \cite{PaVOys:K-CMiso}). In the second appendix we describe 
 differential
structures and differentiable bimodules associated to unital faithfully flat
pre-torsors.

\section{Preliminaries. Corings, bialgebroids, $\times_A$-Hopf algebras}
\label{sec.prelims}

Throughout the paper, all algebras are over a commutative associative ring $k$
with  
a unit. Categories of right (resp.\ left) modules of an algebra $A$ are
denoted by  
$\cM_A$ (resp.\ ${}_A\cM$).

Given an algebra $A$, a coalgebra $\cC$ in the monoidal category of 
$A$-$A$ bimodules is called an {\em $A$-coring}. The coproduct and 
counit in $\cC$ are denoted by $\DC:\cC\to \cC\sstac{A}\cC$ and 
$\eC: \cC\to A$, respectively.  
If a right $A$-module $T$ is a right $\cC$-comodule, then
($A$-linear, coassociative and counital) coaction is denoted by $\roT$. A left
$\cC$-coaction 
in a left $\cC$-comodule $T$ is denoted by $\Tro$. 
The category of right (resp.\ left) $\cC$-comodules is denoted by $\cM^\cC$
(resp.\ ${}^\cC\!\cM$).

Given an $A$-coring $\cC$ and
a $B$-coring $\cD$, a $B$-$A$ bimodule $T$ is called a {\em $\cD$-$\cC$
bicomodule} if $T$ is a right $\cC$-comodule and left $\cD$-comodule with 
$B$-$A$ bilinear coactions
which commute in the sense that
$$
(\Tro\stac{A}\cC)\circ\roT = (\cD\stac{B}\roT)\circ\Tro.
$$
The category of $\cD$-$\cC$ bicomodules is denoted by ${}^\cD\!\cM^\cC$.
In particular, $B$ is a trivial $B$-coring with the coproduct and counit given
by the identity map. Thus to say that $T$ is a $B$-$\cC$ bicomodule is the same
as to say that $T$ is a $B$-$A$ bimodule with left $B$-linear right
$\cC$-coaction.  

For a right $\cC$-comodule $T$ and a left $\cC$-comodule $M$, the cotensor 
product is defined as the equaliser of $\roT\sstac{A}M, T\sstac{A}\Mro:
T\sstac{A}M \to T\sstac{A}\cC\sstac{A}M$ and is denoted by $T\Box_\cC M$.
For a $\cD$-$\cC$ bicomodule $T$ and a left $\cC$-comodule $M$, the cotensor
product $T\Box_\cC M$ is a left $\cD$-comodule provided that $T\Box_\cC M$
is a $\cD\sstac B \cD$-pure equaliser in ${}_B \cM$
(cf. \cite[22.3, erratum]{BrzWis:cor}). If this purity condition holds for
every left $\cC$-comodule $M$, then $T\Box_\cC\, \bullet$ defines a functor 
${}^\cC\cM\to {}^\cD\cM$. The cotensor functor induced by any $\cD$-$\cC$
bicomodule $T$ exists in particular if $\cD$ is a flat right $B$-module.

For more information on corings, the reader is referred to \cite{BrzWis:cor}. 

An algebra $T$ in the monoidal category of $A$-$A$ bimodules is called an 
{\em $A$-ring}. An $A$-ring $T$ is equivalent to a $k$-algebra $T$ and a
$k$-algebra map $\eta_T: A\to T$ (which serves as the {\em unit map} for the 
$A$-ring $T$). 
The unit element in the $k$-algebra $T$ is denoted by $1_T$. The
product in $T$ both as a map $T\sstac{k}T\to T$ and $T\sstac{A}T\to T$
is denoted by $\mu_T$. 

A triple $(T,\cC,\psi)$, where $T$ is an $A$-ring,
$\cC$ is an $A$-coring and $\psi: \cC\sstac{A}T\to T\sstac{A}\cC$ is an $A$-$A$
bilinear map such that
\begin{eqnarray*}
\psi\circ (\cC\stac{A} \mu_T)&=&(\mu_T\stac{A}
\cC)\circ(T\stac{A}\psi)\circ (\psi \stac{A} T),\qquad
\psi\circ(\cC\stac{A} \eta_T)=\eta_T\stac{A} \cC,\\
(T\stac{A} \Delta_\cC)\circ \psi &=& (\psi\stac{A} \cC) \circ
(\cC\stac{A} \psi)\circ (\Delta_\cC\stac{A} T), \qquad
(T\stac{A} \eC)\circ \psi =\eC \stac{A} T,
\end{eqnarray*}
is called a {\em right entwining structure over $A$}. A {\em right entwined
 module} is a a right $\cC$-comodule
and right $T$-module $M$, with multiplication $\varrho_M: M\sstac{A}T\to M$, 
such that
$$
\varrho^M\circ \varrho_M=(\varrho_M\stac{A}\cC)\circ (M\stac{A}\psi)\circ
(\varrho^M\stac{A} T).
$$
The morphisms of entwined modules are right $A$-linear, right
$\cC$-colinear maps. If $T$ itself is an entwined module with multiplication
$\mu_T$ and coaction $\roT$, then, for any entwined module $M$, the {\em
coinvariants} are defined by
$$
M^{co\cC} := \{m\in M \; |\; \varrho^M(m) = m\roT(1_T)\}.
$$
In particular, to any $A$-coring $\cC$ one can associate a 
{\em trivial right entwining structure} $(A,\cC,\cC)$ (i.e.\ the entwining map
$\psi$ is the identity map on $\cC$). Entwined modules in this case
are simply right comodules of $\cC$. Thus to say that $A$ is an entwined 
module in this case is the same as to say that $\cC$ has a group-like 
element $g$, and $\roA(1_A) = 1_A\sstac{A}g$. Consequently the coinvariants of 
a right $\cC$-comodule 
$N$ coincide with $\{n\in N\; |\; \varrho^N(n) = n\sstac{A}g\}$.

Every right entwining structure $(T,\cC,\psi)$ over $A$ gives rise to a
$T$-coring 
$T\sstac{A}\cC$. Entwined modules can be identified with right comodules of
this $T$-coring. 

Left entwining structures are defined symmetrically. In
particular if $\psi$ in a right entwining structure $(T,\cC,\psi)$  is
bijective, then $(T,\cC,\psi^{-1})$ is a left entwining structure. 

Let $\cC$ be an $A$-coring and $T$ be an $A$-ring which is also a right 
$\cC$-comodule (with right $A$-module structure given by the unit map). Define
$B= \{b\in T \; |\; \forall t\in T\; \roT(bt) = b\roT(t)\}$. We say that $T$
is a {\em right $\cC$-Galois extension of $B$} if the canonical Galois map
$$
\can : T\stac{B}T\to T\stac{A}\cC, \qquad t\stac{B}t'\mapsto t\roT(t'),
$$
is bijective. The restricted inverse of $\can$,
$
\chi : \cC\to T\sstac{B}T$, $c\mapsto \can^{-1}(1_T\sstac{A}c),
$
is called the {\em translation map}. 
 If $T$ is a left $\cC$-comodule, one
defines 
a {\em left $\cC$-Galois extension} in analogous way, using the 
left coaction to define the {\em left} Galois map. When dealing with both
right and left $\cC$-Galois extensions we write $\can_\cC$ for the right
Galois map and ${}_\cC\can$ for the left Galois map.

Every right (resp.\ left) $\cC$-Galois extension $B\subseteq T$ gives rise to
a right (resp.\ left) entwining structure $(T,\cC,\psi)$ with 
$\psi: c\sstac A t\mapsto \can_\cC\big(\can_\cC^{-1}(1_T\sstac A c)\ t\big)$
(resp. $\psi: t\sstac A c \mapsto {}_\cC \can\big(t\, {}_\cC\can^{-1}(c\sstac A
1_T)\big)$), for which $T$ is an entwined module. One then easily shows that
$B=T^{co\cC}$.  

An $A$-coring $\cC$ is a  {\em Galois coring} if $A$ itself is a $\cC$-Galois
extension of the coinvariant subalgebra $\{b\in A \; |\; \forall a\in A\;
\roA(ba) = b\roA(a)\}$ (which is the commutant in the $A$-$A$ bimodule $\cC$
of the grouplike element $g$ determining the $\cC$-coaction in $A$).
The $T$-coring $T\sstac A \cC$, arising from an entwining structure determined
by a $\cC$-Galois extension $B\subseteq T$, is a Galois coring (with grouplike
element $\roT(1_T)$). The coinvariants of $T$ with respect to $\cC$ and
$T\sstac A \cC$ are the same.

For a $T$-coring ${\mathcal E}$ with a grouplike element, there exist adjoint
functors 
\begin{equation}\label{eq:gadj}
(\bullet)^{co{\mathcal E}}:\cM^{\mathcal E}\to \cM_B\quad \textrm{and}\quad
  \bullet \stac B T: \cM_B \to \cM^{\mathcal E},
\end{equation}
(cf. \cite[28.8]{BrzWis:cor}) where $B\colon = T^{co{\mathcal E}}$ is the
coinvariant subalgebra. By the Galois Coring Structure Theorem \cite[28.19 (2)
  (a) $\Rightarrow$ (c)]{BrzWis:cor}, the functors \eqref{eq:gadj} are inverse
equivalences provided that ${\mathcal E}$ is a Galois coring and $T$ is a
faithfully flat left $B$-module.

\begin{lemma}
\label{lem:purecoC}
Let ${\mathcal E}$ be a Galois coring over an algebra $T$ and let $A$ be an
algebra.  Assume that $T$ is a faithfully flat left module for its ${\mathcal
  E}$-coinvariant subalgebra  $B$. For any $A$-${\mathcal E}$ bicomodule $M$, 
$M^{co{\mathcal E}}$ is a pure equaliser in ${}_A\cM$.
\end{lemma}

\begin{proof}
For any right $A$-module $N$, $N\sstac A M$ inherits a right ${\mathcal
  E}$-comodule structure of $M$. Applying twice the Galois Coring 
  Structure Theorem \cite[28.19 (2) $(a)\Rightarrow (c)$]{BrzWis:cor}, one
  concludes that
$$
(N\stac A M)^{co{\mathcal E}} \stac B T\cong N\stac A M \cong N\stac A
M^{co{\mathcal E}}\stac B T.
$$
Hence the claim follows by the faithful flatness of $T$ as a left $B$-module.
\end{proof}

By an analogy to  rings, a $B$-coring $\cD$ is termed 
a {\em right extension} of an $A$-coring $\cC$ if the forgetful
functor $\cM^\cC\to \cM_k$ factorises through a $k$-linear functor $\cM^\cC\to
\cM^\cD$ and the forgetful functor $\cM^\cD\to \cM_k$. In 
\cite[Theorem 2.6]{Brz:corext} this definition was shown to be equivalent to
the existence of a right $\cD$-coaction on $\cC$, which is left colinear
with respect to the left regular comodule structure of $\cC$. Left extensions
of corings are defined symmetrically, in terms of the categories of left
comodules. 


Extensions of corings over the same base algebra are related to coring
morphisms as follows.\footnote{Lemma~\ref{lem:correct} and new formulation of
  Lemma~\ref{lemma.unique} replace incorrect statements on page 548 and in
  Lemma~3.7 of the published version of this paper:  [G.\ B\"ohm \& T.\
    Brzezi\'nski, {\em Pre-torsors and  equivalences}, J. Algebra 317 (2007)
    544--580].} 
\begin{lemma}\label{lem:correct}
For two $A$-corings ${\mathcal C}$ and ${\widetilde {\mathcal C}}$, denote the
 forgetful functors by $F:{\mathcal M}^{\mathcal C} \to {\mathcal M}_A$ 
and ${\widetilde F}:{\mathcal M}^{\widetilde{\mathcal C}} \to {\mathcal
 M}_A$. Then the following assertions are equivalent.
\begin{rlist}
\item There is a $k$-linear functor $U:{\mathcal M}^{\mathcal C} \to 
  {\mathcal M}^{\widetilde {\mathcal C}}$ such that $F= {\widetilde F} \circ
  U$.  
\item Considering the $A$-bimodule ${\mathcal C}$ as a left ${\mathcal
  C}$-comodule via the coproduct, there is a right ${\widetilde{\mathcal
  C}}$-coaction ${\widetilde \varrho}:{\mathcal C} \to {\mathcal C} \otimes_A
  {\widetilde {\mathcal C}}$, making ${\mathcal C}$ a ${\mathcal
  C}$-${\widetilde  {\mathcal C}}$ bicomodule.  
\item There is a homomorphism of $A$-corings $\kappa:{\mathcal C} \to
  {\widetilde {\mathcal C}}$.
\end{rlist}
\end{lemma}
\begin{proof}
(i)$\Rightarrow$ (ii) By property (i), there is a right ${\widetilde{\mathcal
    C}}$-coaction ${\widetilde \varrho}$ on the right $A$-module ${\mathcal
    C}$. Since under assumption (i) ${\widetilde {\mathcal C}}$ is a right
    extension of ${\mathcal C}$, ${\widetilde \varrho}$ is a left ${\mathcal
    C}$-comodule map by \cite[Theorem 2.6]{Brz:corext}. 

(ii) $\Rightarrow$ (iii) The map $\kappa$ is constructed
as $\kappa:= (\epsilon_{\mathcal C}\otimes_A {\widetilde {\mathcal C}})\circ
{\widetilde \varrho}$. 

(iii) $\Rightarrow$ (i) The functor $U$ is given by the {\em corestriction
  functor} along $\kappa$, cf. \cite[22.11]{BrzWis:cor}.
\end{proof}

A symmetrical statement holds for the categories of left (co)modules. 
Lemma \ref{lem:correct} (ii) asserts that the $A$-module structures of
${\mathcal C}$, as an $A$-coring on one hand and as a ${\tilde {\mathcal 
C}}$-comodule on the other, are the same. Note that this property allows 
the construction in the proof of the implication (ii)
$\Rightarrow$ (iii) in Lemma \ref{lem:correct} to yield a well
defined and right $A$-linear map $\kappa$.

A right bialgebroid  over $A$ \cite{Tak:crR},\cite{Lu:hgd} is a quintuple 
$\cC=(\cC, s,t,\DC,\eC)$. Here $s,t: A\to \cC$ are $k$-linear maps,
$\cC$ is an $A\sstac{k}A^{op}$-ring with
the unit map $\mu_\cC\circ(s\sstac{k}t)$ (thus, in particular, $s$ is an
algebra and $t$ is an anti-algebra morphism), and $(\cC,\DC,\eC)$ is an
$A$-coring. The bimodule structure of
this $A$-coring is given by
$$
aca'=cs(a')t(a),\qquad \textrm{for all }a,a'\in A,\ c\in \cC.
$$
The range of the coproduct is required to be in the {\em Takeuchi product}
$$ \cC\times_A \cC :=
\{\ \sum_i c_i\stac{A} c'_i\in \cC\stac{A} \cC\; |\;   \forall a\in A,\; 
\sum_i s(a)c_i \stac{A} c'_i =\sum_i c_i\stac{A} t(a)c'_i \
 \},
$$
which is an algebra by factorwise multiplication, and the (corestriction of
the) coproduct is
required to be an algebra map. The  counit satisfies the following conditions
$$
\eC(1_\cC) = 1_A, \quad \eC\left(s(\eC(c))c'\right)=
\eC\left(t(\eC(c))c'\right)= 
\eC(cc'),
$$
for all $c,c'\in \cC$. The map $s$ is called a {\em source map}, and $t$ is
known as a {\em target map}.
Left bialgebroids are defined analogously in terms of multiplications by the
source and target maps on the
left. For more details we refer to \cite{KadSzl:D2bgd}.

The category of right comodules of a right $A$-bialgebroid $(\cC, s,t,\DC,\eC)$
is a monoidal category with a strict monoidal functor to the category of
$A$-$A$  bimodules  \cite[Proposition 5.6]{Sch:Bianc}. Explicitly, 
take a right comodule $M$ of $\cC$ (i.e. of the $A$-coring $(\cC,\DC,\eC)$),
and define the left $A$-multiplication on $M$ in terms of the coaction 
$\roM: M\to M\sstac{A}\cC$,  $m\mapsto m\su 0\sstac{A} m\su 1$, 
$$
am=m\su 0 \eC\big(t(a)m\su 1\big) =  m\su 0 \eC\big(s(a)m\su 1\big),
\qquad \textrm{for all }a\in A,\ m\in M.
$$
This is a unique left multiplication which guarantees that the
image of the coaction $\roM$ is a subset of the Takeuchi product
$$
M\times_A \cC :=
\{\ \sum_i m_i\stac{A} c_i\in M\stac{A} \cC\; |\; \forall\ a\in A,\;
\sum_i am_i\stac{A} c_i= \sum_i m_i\stac{A} t(a)c_i
 \}.
$$
In addition it equips $M$ with an $A$-$A$ bimodule structure such that every
$\cC$-colinear map becomes $A$-$A$ bilinear.
For any right $\cC$-comodules $M$ and $N$, the right $A$-multiplication and the
right $\cC$-coaction on the
tensor product $M\sstac{A}N$ are defined by
$$
(m\stac{A}n)\cdot a = m\stac{A} na,\quad (m\stac{A} n)\su 0 \stac{A} 
(m\stac{A}
  n)\su 1 =
(m\su 0\stac {A} n\su 0)\stac{A} m\su 1 n\su 1,
$$
for all $a\in A$, $m\sstac{A} n\in M\sstac{A} N$. The monoidal unit is $A$,
with the regular module structure and coaction given by the source map. In a
symmetric way, also 
the left comodules of a right $A$-bialgebroid $(\cC, s,t,\DC,\eC)$ form a
monoidal category with respect to the $A^{op}$-module tensor product. 
Any left $\cC$-comodule $M$, with coaction ${}^M\varrho:M\to
\cC\sstac A M$, $m\mapsto m\su {-1}\sstac A m\su 0$, is equipped with a right
$A$-module structure 
$$
ma = \varepsilon_\cC(s(a) m\su {-1})m\su 0= \varepsilon_\cC(t(a) m\su
{-1})m\su 0,\quad \textrm{for all}\quad a\in A,\ m\in M.
$$
With this definition, the range of ${}^M\!\varrho$ is in the Takeuchi product
$$ 
\cC\times_A M:=\{\ \sum_i c_i\stac A m_i \in \cC\stac A M\ |\ 
\forall \ a\in A,\ \sum_i
c_i\stac A m_i a = \sum_i s(a) c_i \stac A m_i\ \},
$$ 
and $M$ is an $A$-$A$ (equivalently, $A^{op}$-$A^{op}$) bimodule such that any
comodule map is bilinear. For any two left $\cC$-comodules $M$ and $N$, with
coactions $m \mapsto m\su {-1}\sstac A m\su 0$ and $n\mapsto n\su {-1}\sstac A
n\su 0$, respectively, the left $A$-action and left $\cC$-coaction are 
$$
a(m\stac {A^{op}} n)=m\stac {A^{op}} an, \qquad 
(m\stac {A^{op}} n)\su{-1}\stac A (m\stac {A^{op}} n)\su 0 =
m\su{-1} n\su{-1}\stac A (m\su 0 \stac{A^{op}} n\su 0),
$$
for $a\in A$ and $m\sstac {A^{op}} n\in M\sstac {A^{op}} N$. The monoidal unit
is $A^{op}$, with the left regular $A$-module structure and coaction given by
the target map.  

A right {\em comodule algebra} $T$ for a right bialgebroid $(\cC, s,t,\DC,\eC)$
over $A$ is an algebra in the monoidal category of right comodules for $(\cC,
s,t,\DC,\eC)$. By the existence of a strict monoidal forgetful functor from
this comodule category to the category of $A$-$A$ bimodules, $T$ is in
particular an $A$-ring. The $A$-ring $T$ and the $A$-coring $(\cC,\DC,\eC)$
are canonically entwined. A right {\em Galois extension} by a  right
bialgebroid $(\cC, s,t,\DC,\eC)$ over $A$ is a right comodule algebra $T$, 
which is a Galois extension of its coinvariant subalgebra by the $A$-coring
$(\cC,\DC,\eC)$. The range of the translation map in this case is contained 
in 
the center of the
$B$-$B$ bimodule $T\sstac B T$.

Note that the definition of the left $A$-multiplication on a right
$\cC$-comodule 
$T$ implies that, for any left $\cC$-comodule $M$, the cotensor product 
$T\Box_\cC M$ is contained in the $A$-centraliser of $T\sstac{A}M$. 
That is, denote  by $t\mapsto t\su 0 \sstac A t\su 1$ and $m\mapsto
m\su{-1}\sstac A m\su 0$ the right coaction in $T$ and the left coaction
in $M$, respectively. For $\sum_i t_i\sstac A m_i\in T\Box_\cC M$ and $a\in A$,
\begin{eqnarray*}
a(\sum_i t_i\stac A m_i)
&=&\sum_i t_i \su 0 \varepsilon_\cC(s(a) t_i\su 1)\stac A m_i\\
&=&\sum_i t_i \stac A \varepsilon_\cC(s(a) m_i \su {-1}) m_i \su 0
=(\sum_i t_i\stac A m_i)a.
\end{eqnarray*}
The second equality follows by the definition of the cotensor product as
an equaliser.

Given a right bialgebroid $(\cC, s,t,\DC,\eC)$ over $A$, one can define the
Galois map corresponding to the right regular $\cC$-comodule, with
coinvariants $t(A)\cong A^{op}$, 
\begin{equation}\label{eq:can_x_A}
\theta: \cC\stac {A^{op}} \cC\to \cC\stac
A \cC, \qquad c\stac A c'\mapsto c\Delta_\cC(c').
\end{equation}
The map $\theta$ satisfies the pentagon identity
\begin{equation}\label{eq:pent}
(\theta\stac A \cC)\circ (\cC \stac {A^{op}} \theta)=(\cC \stac A \theta)\circ
\theta_{13}\circ (\theta\stac {A^{op}} \cC),
\end{equation}
where $\theta_{13}:(\cC_\bullet \sstac{A}\cC)\sstac{A^{op}}\cC\to 
\cC\sstac{A}(\cC\sstac{A^{op}}{}_\bullet\cC)$ is the map defined
as the nontrivial action of $\theta$ on the first and the third factors.
 Following \cite{Sch:dua}, a right bialgebroid $(\cC, s,t,\DC,\eC)$ over $A$
 is called  
a {\em $\times_A$-Hopf algebra} iff the Galois map $\theta$ is bijective.
Similarly, for left $A$-bialgebroids, one considers the left Galois map
 $c\sstac{A^{op}}c'\mapsto \DC(c)c'$ and defines a left $\times_A$-Hopf 
algebra by requiring this map be bijective. The translation map
$\theta^{-1}(1_\cC\sstac A \, \bullet)$, corresponding to \eqref{eq:can_x_A},
is an algebra map from $\cC$ to the center of the $A^{op}$-$A^{op}$ bimodule
 $\cC\sstac {A^{op}} \cC$, where all $A^{op}$-module structures are given by
 the target map and the algebra structure is inherited from $\cC^{op}\sstac k
 \cC$. Furthermore,
\begin{equation}\label{eq:thetaA}
\theta^{-1}(1_\cC\stac A t(a))=s(a) \stac {A^{op}}  1_\cC\quad
\textrm{and}\quad 
\theta^{-1}(1_\cC\stac A s(a))=1_\cC  \stac {A^{op}} s(a),
\end{equation}
for all $a\in A$. 

Let $(\cC, s,t,\DC,\eC)$ be a right bialgebroid over $A$ which admits a  
Galois extension $B\subseteq T$ such that $T$ is a faithfully flat right 
$A$-module. By \cite[Lemma 4.1.21]{Hobst:PhD}, the existence of such an
extension implies that $(\cC, s,t,\DC,\eC)$ is a $\times_A$-Hopf 
algebra.

\section{Pre-torsors, coring-Galois extensions and 
bi-Galois objects}\label{sec.pre-tors}
One of the most striking properties of  $A$-$B$ torsors is the observation
made in  
 \cite[Theorem~5.2.10]{Hobst:PhD} that, in the faithfully flat case,  
they correspond  
to Galois extensions by bialgebroids. Guided by the non-commutative geometry 
experience, where the notion of a Hopf-Galois extension is not flexible enough
to  
describe examples of quantum principal bundles, one can also envisage 
that the notion of an $A$-$B$ torsor might be too strict to deal with bundles
over 
quantum groupoids (cf.\ Example~\ref{ex.homog}). 
In this section we introduce the notion of an 
{\em $A$-$B$ pre-torsor}. In the faithfully flat case we show that this notion
is equivalent 
to a certain class of coring-Galois extensions. 

\begin{definition}\label{def.pre.tor}
Let $\alpha: A\to T$ and $\beta: B\to T$  be 
$k$-algebra maps. View $T$ 
as an $A$-$B$ bimodule and a $B$-$A$ bimodule via the maps $\alpha$ and 
$\beta$. We say that $T$ is an {\em $A$-$B$ pre-torsor} if there exists a 
$B$-$A$ bimodule map
$$
\tau : T\to T\stac{A}T\stac{B}T, \qquad t\mapsto t\sut 1\stac{A} t\sut 2
\stac{B}  t\sut 3,
$$
(summation understood) such that
\begin{blist}
\item $(\mu_T\stac{B}T)\circ \tau = \beta\stac{B}T$;
\item $(T\stac{A}\mu_T)\circ \tau = T\stac{A} \alpha$;
\item $ (\tau \stac{A}T\stac{B}T)\circ \tau = (T\stac{A}T\stac{B}\tau)\circ
  \tau$, 
\end{blist}
where $\mu_T$ denotes the quotients of the product in $T$ to appropriate
tensor products (over $A$ or $B$).

An $A$-$B$ pre-torsor is said to be {\em faithfully flat}, if it is faithfully
flat as a right $A$-module and left $B$-module. (Note that in this case the
unit maps $\alpha$ and $\beta$ are injective.) 

\end{definition}
Following the observation made in \cite[Remark~2.4]{BohBrz:str} we propose the 
following
\begin{definition}\label{def.biGal}
Given an $A$-coring $\cC$ and a $B$-coring $\cD$, an algebra $T$ is called 
a {\em $\cC$-$\cD$ bi-Galois object} if $T$ is a $\cD$-$\cC$ bicomodule, and 
$T$ is a right $\cC$-Galois extension of $B$ and a left $\cD$-Galois extension 
of $A$. A $\cC$-$\cD$ bi-Galois object $T$ is said to be {\em faithfully
  flat}, if $T$ 
is faithfully flat as a left $B$-module and a right $A$-module.
\end{definition}

Since corings can be understood as the algebraic structure underlying
quantum (Lie) groupoids, the $\cC$-Galois extension $B\subseteq T$
can be understood as  the dual to a free action of a groupoid on a manifold.
With this interpretation in mind a (faithfully flat)  $\cC$-$\cD$ bi-Galois
object $T$ 
can be seen as a non-commutative version of a groupoid principal bundle over
a groupoid (cf.\ \cite[Section~5.7]{MoeMrc:fol}).

The main result of this section, which includes the pre-torsor version of
\cite[Theorem~5.2.10]{Hobst:PhD}, is contained in Theorem \ref{thm.main.1}
below.  

\begin{remark}\label{rem:rel_to_Hobst}
The basic difference between \cite[Theorem~5.2.10]{Hobst:PhD} and 
Theorem~\ref{thm.main.1} is that the former one is formulated for $A$-$B$
torsors, 
while the latter one deals with pre-torsors, i.e.\ without multiplicative
structure of the torsor map. Furthermore, 
\cite[Theorem~5.2.10]{Hobst:PhD} deals with the analogues of the implications
(i)$\Rightarrow$ (ii) and (i)$\Rightarrow$ (iii) in Theorem~\ref{thm.main.1} 
not with equivalences. Also, in \cite[Theorem~5.2.10]{Hobst:PhD} 
a torsor $T$ is
assumed to be faithfully flat both as a left and right module, for both base
algebras $A$ and $B$. In contrast, we assume only the faithful flatness of a
pre-torsor $T$ as a right $A$-module and as a left 
$B$-module, cf.\ Definition~\ref{def.pre.tor}. 
Thanks to Lemma~\ref{lemma.pure.eq}, these two faithful
flatness assumptions are enough to conclude the claim. 
\end{remark}

\begin{theorem} \label{thm.main.1}
There is a bijective correspondence between the following sets of data:
\begin{rlist}
\item faithfully flat $A$-$B$ pre-torsors $T$;
\item $A$-corings $\cC$ and left faithfully flat right $\cC$-Galois extensions
$B\subseteq T$, such that $T$ is a right faithfully flat $A$-ring;
\item $B$-corings $\cD$ and right faithfully flat left $\cD$-Galois extensions 
$A\subseteq T$, such that $T$ is a left faithfully flat $B$-ring.
\end{rlist}
Furthermore, a faithfully flat  $A$-$B$ pre-torsor $T$ is a faithfully
flat $\cC$-$\cD$ bi-Galois object with $\cC$ and $\cD$ as in parts (ii) and
(iii). 
\end{theorem}

\begin{proof} 
We will prove the equivalence of statements (i) and (ii). The equivalence of 
(i) with (iii) will
then follow by the symmetric nature of the notion of an $A$-$B$ pre-torsor.
The implication (ii)~$\Ra$~(i) is a  consequence of the following

\begin{lemma}\label{lemma.tor1}
Let $\cC$ be an $A$-coring and $T$ an $A$-ring. 
If $B\subseteq T$ is a right $\cC$-Galois extension, then
$T$ is an $A$-$B$ pre-torsor with the structure map 
$$
\tau := (T\stac{A}\chi)\circ\varrho^T: T\to T\stac{A}T\stac{B}T,
$$
where $\varrho^T: T\to T\sstac{A}\cC$ is the coaction and $\chi: \cC\to
T\sstac{B} T$ is the translation map. 
\end{lemma}
\begin{proof}
Since $\varrho^T$ is a $B$-$A$ bimodule map and $\chi$ is right $A$-linear,
the map $\tau$ is $B$-$A$ bilinear.  Note that $(\mu_T\sstac{B}T)\circ\tau =
(\can)^{-1}\circ 
\varrho^T$, and $\varrho^T = \can\circ(\beta\sstac{B}T)$, where $\beta$ is
the  
obvious inclusion $B\subseteq T$. This implies that $\tau$ satisfies condition
(a) in  
Definition~\ref{def.pre.tor}. The condition (b) follows by noting that
$\alpha\circ \eC =  \mu_T\circ\chi$,
where $\alpha$ is the unit map $A\to T$. 
A simple calculation which uses the coassociativity of $\varrho^T$ and
right $\cC$-colinearity of $\chi$ confirms that $\tau$ satisfies property (c).
\end{proof}

Now assume that $T$ is a faithfully flat $A$-$B$ pre-torsor with the
structure map $\tau$. The proof of the converse implication (i)~$\Ra$~(ii)
starts with the following 

\begin{lemma}\label{lemma.pure.eq}
The equaliser:
\begin{equation}\label{eq:C_eq}
\xymatrix{
  \cC\ar[rr] &&T\stac{B}T 
  \ar@<.6ex>[rrrr]^{(\mu_T\stac{A}T\stac{B}T)\circ (T\stac{B}\tau)}
   \ar@<-.6ex>[rrrr]_{\alpha\stac{A}T\stac{B}T}  &&&&
    T\stac{A}T\stac{B}T} 
\end{equation}
is pure in $\M_A$.
\end{lemma}
\begin{proof}
\footnote{The referee pointed out an alternative proof of Lemma 3.6. It is
  based on the observation that the equaliser, obtained by tensoring (3.1)
  with the faithfully flat right A-module T on the left, is contractible hence 
  cosplit.}
Let $\omega = (\mu_T\sstac{A}T\sstac{B}T)\circ (T\sstac{B}\tau) - 
\alpha\sstac{A}T\sstac{B}T$. For any left $A$-module $N$, define a map
$
\varphi_N: T\sstac{A}\ker (\omega\sstac{A}N)\to T\sstac{B} T\sstac{A}N,
$
as the restriction of $\mu_T\sstac{B}T\sstac{A}N$.  Consider also the map
$\widetilde{\varphi}_N := 
(\mu_T\sstac{A}T\sstac{B}T\sstac{A}N)\circ(T\sstac{B}\tau\sstac{A}N).
$
Using properties Definition~\ref{def.pre.tor}~(b) and (c) one easily finds that
$ (T\sstac{A}\omega\sstac{A}N)\circ \widetilde{\varphi}_N =0$, hence, in view
of the  
flatness of $T$ as a right $A$-module, $\im \widetilde{\varphi}_N \subseteq
T\sstac{A}\ker (\omega\sstac{A}N)$. The equality $\varphi_N\circ
\widetilde{\varphi}_N 
= T\sstac{B}T\sstac{A}N$ follows by  Definition~\ref{def.pre.tor}~(a), while 
the definition of $\omega$ implies that the composition
$\widetilde{\varphi}_N\circ  
\varphi_N$ is the identity too. Thus $\varphi_N$ is a left $T$-module 
isomorphism (which is $T$-$A$-linear if $N$ is an $A$-bimodule) 
and hence the composition
$$
\xymatrix{
T\stac{A}\ker\omega\stac{A}N \ar[rr]^{\varphi_A\stac{A}N} &&
T\stac{B}T\stac{A}N 
\ar[rr]^{\varphi_N^{-1}} &&  T\stac{A} \ker (\omega\stac{A}N) 
 ,}
$$
which is simply the obvious map, 
is an isomorphism too. By the faithful flatness of $T$, 
this yields the 
required isomorphism
$\ker\omega\sstac{A}N \cong \ker(\omega\sstac{A} N) $.
\end{proof}
 
 The equaliser $\cC$ in Lemma~\ref{lemma.pure.eq} is given in terms of
 $A$-bimodule 
 maps, hence it is an $A$-bimodule with the $A$-multiplications given through
 $\alpha$,
 $$
 a\left(\sum_i t_i\stac{B}u_i\right) a' := \sum_i
 \alpha(a)t_i\stac{B}u_i\alpha(a'). 
 $$
 Define an $A$-bimodule map
 $$
 \DC :\cC \to T\stac{B}T\stac{A}T\stac{B}T, \qquad \DC = T\stac{B}\tau.
 $$ 
 Since $\tau = \widetilde{\varphi}_A \circ (\beta\sstac{B}T)$ (cf.\ the proof
 of 
 Lemma~\ref{lemma.pure.eq}), the range of $\tau$ is contained in
 $T\sstac{A}\cC$,  
 hence $\DC(\cC)\subseteq T\sstac{B}T\sstac{A}\cC$. 
 Using the definition of $\cC$ as the kernel of the map $\omega$ (cf.\ the
 proof of 
 Lemma~\ref{lemma.pure.eq}) and property (c) in Definition~\ref{def.pre.tor},
 one 
 immediately finds that $\im \DC \subseteq \ker(\omega\sstac{A}\cC)$, hence
 $\DC(\cC)\subseteq \cC\sstac{A}\cC$ by  Lemma~\ref{lemma.pure.eq}. 
 
 The map $\DC$ is coassociative by Definition~\ref{def.pre.tor}~(c). Note that 
 property (b) in Definition~\ref{def.pre.tor} implies that 
 $(T\sstac{A}\mu_T)\circ \omega = \mu_T\sstac{A}\alpha -
 \alpha\sstac{A}\mu_T$.   This means that, for all $c\in \cC$,
 $\mu_T(c)\sstac{A}1_T=1_T\sstac{A}\mu_T(c)$.  
 In view of
 the faithful flatness of the right $A$-module $T$, this implies that
there is an $A$-$A$ bimodule map
 $$
 \eC: \cC\to A, \qquad \sum_i t_i\stac{B}u_i\mapsto \sum_i t_iu_i.
 $$
The map $\eC$ is a counit for $\DC$ by properties (a) and (b) in 
Definition~\ref{def.pre.tor}.

The torsor map
$\tau: T\to T\sstac{A}\cC$ is a right $\cC$-coaction of $\cC$ on $T$, which is
coassociative by property (c) in Definition~\ref{def.pre.tor}, and counital by
Definition~\ref{def.pre.tor}~(b). By definition, $\tau$ is left
$B$-linear. Conversely, if  
$b\in T$ is such that, for all $t\in T$, $\tau(bt) = b\tau(t)$, then property
(a) in  
Definition~\ref{def.pre.tor} implies that $1_T\sstac{B}bt = b\sstac{B}t$,
hence  
$b\in \im\,\beta$ by the faithful flatness of $T$ as a left $B$-module. 
Finally note that the canonical map $T\sstac B T\to T\sstac A \cC$ is the same
as the bijective map $\widetilde{\varphi}_A$ 
constructed in the proof of Lemma~\ref{lemma.pure.eq}. Therefore, $B\subseteq
T$ is a   
right $\cC$-Galois extension as required. 

Thus we have established a correspondence between faithfully flat $A$-$B$
pre-torsors and faithfully flat right $\cC$-Galois extensions. There is a
symmetric correspondence between faithfully flat $A$-$B$ pre-torsors $T$ and 
faithfully flat left $\cD$-Galois extensions $A\subseteq T$, for the
$B$-coring $\cD=\mathrm{ker}((T\sstac A T\sstac B \mu_T)\circ (\tau\sstac A
T)-T\sstac A T \sstac B \beta)$. Both the right $\cC$-coaction  
and the left $\cD$-coaction in $T$ are given by the torsor map $\tau$, hence
$T$ is a $\cD$-$\cC$ bicomodule by property (c) in
Definition~\ref{def.pre.tor}.  

It remains to prove that the established correspondence between faithfully
flat $A$-$B$ pre-torsors and faithfully flat coring-Galois extensions
is bijective. This is a consequence of the following
\begin{lemma}\label{lemma.unique}
Let $T$ be a faithfully flat $A$-$B$ pre-torsor and ${\mathcal C}$
and ${\mathcal D}$ be the associated $A$- and $B$-corings, respectively. Let 
${\widetilde {\mathcal C}}$ be an $A$-coring for which $T$  is a ${\mathcal
  D}$-${\widetilde {\mathcal C}}$ bicomodule. Then the following hold: 
\begin{zlist}
\item  If the right $A$-action of $T$ as a right ${\widetilde {\mathcal
  C}}$-comodule is determined by its $A$-ring structure, then there is a
  homomorphism of $A$-corings $\kappa:{\mathcal C} \to  {\widetilde {\mathcal
  C}}$.  
\item $B\subseteq T$ is a right $\widetilde\cC$-Galois extension if and only
  if the coring homomorphism $\kappa$ in part (1) is an isomorphism. 
\end{zlist}
\end{lemma}
\begin{proof}
(1) In view of Lemma \ref{lem:correct} (ii) $\Rightarrow$(iii), we need to
    construct a right 
    $\widetilde \cC$-coaction on the right $A$-module $\cC$ which is 
    left colinear with respect to the regular comodule structure of
    $\cC$. We claim that such a coaction is given in terms of the
    $\widetilde\cC$ coaction in $T$, $t \mapsto t^{\widetilde{[0]}}\sstac A
    t^{\widetilde{[1]}}$, as 
\begin{equation}
\label{eq:tC}
\cC\to \cC\stac A {\widetilde \cC}, \qquad \sum_i u_i\stac B v_i \mapsto 
\sum_i u_i \stac B {v_i}^{\widetilde{[0]}}\stac A {v_i}^{\widetilde{[1]}}.
\end{equation}
Applying the map $\omega\sstac A {\widetilde \cC}$ in Lemma
\ref{lemma.pure.eq} to the range of \eqref{eq:tC}, 
and using the $\cD$-$\widetilde \cC$ bicomodule property of $T$ together with
the characterisation of $\cC$ as the kernel of $\omega$, we conclude that 
the range of the map \eqref{eq:tC} is in $\mathrm{ker}(\omega\sstac A
  {\widetilde \cC})$. Hence 
by Lemma \ref{lemma.pure.eq} it is in $\cC\sstac A {\widetilde \cC}$, as
required. Counitality 
and coassociativity of the coaction \eqref{eq:tC} hold by the counitality
and coassociativity of the $\widetilde \cC$-coaction in $T$. Its left
$\cC$-colinearity follows immediately by the $\cD$-$\widetilde \cC$ bicomodule
property of $T$. Since the right $A$-multiplication in $T$ is determined by
the $A$-ring structure of $T$, the assumptions of
Lemma~\ref{lem:correct}~(ii) are satisfied, and hence  
$$
\kappa:\cC \to {\widetilde{\cC}}, \qquad \sum_i u_i\stac B
 v_i\mapsto \sum_i (u_i v_i^{\widetilde{[0]}})v_i^{\widetilde{[1]}}
$$ 
is a homomorphism of $A$-corings.

(2) Recall that in a ${\widetilde {\mathcal C}}$-Galois
extension $B\subseteq T$, $T$ is an $A$-ring and a ${\widetilde {\mathcal
    C}}$-comodule via the same right $A$-actions, so the homomorphism $\kappa$
in part (1) can be constructed.  It follows by  Definition \ref{def.pre.tor}
(a) that   
$(T\sstac A \kappa)\circ \can_\cC$ is equal to the $\widetilde \cC$-canonical 
map $\can_{\widetilde\cC}$.
Hence if $\kappa$ is an isomorphism of corings, then a $\cC$-Galois extension
$B\subseteq T$ is $\widetilde \cC$-Galois. Conversely,
if $B\subseteq T$ is a right $\widetilde \cC$-Galois extension, i.e.
$\can_{\widetilde\cC}$ is bijective, then the bijectivity of $\kappa$ follows
by the bijectivity of $\can_\cC$ and the faithful flatness of $T$ as a right
$A$-module. 
\end{proof}

This completes the proof of Theorem~\ref{thm.main.1}.
\end{proof}

We conclude this section with a number of examples of $A$-$B$ pre-torsors.
\begin{example}
The simplest examples of Galois extensions by corings are Galois
corings. 
Indeed, for a Galois $T$-coring $\cC$, the base algebra $T$ is a
$\cC$-Galois extension of the coinvariant subalgebra $B:=T^{co\cC}$. 
By Lemma \ref{lemma.tor1}, to this $\cC$-Galois extension there corresponds a
$T$-$B$ pre-torsor $(T,\tau)$, with
$$
\tau: T\to T\stac T T \stac B T\cong T\stac B T,\qquad t\mapsto 1_T\stac B t.
$$
\end{example}

\begin{example}\label{ex.homog}
Examples of Galois extensions by a coring can be obtained by generalising the
construction of quantum homogenous spaces due to Schneider \cite{Schn:PriHomSp}
\cite{MulSch:QHomSpff} from Hopf algebras to $\times_A$-Hopf algebras. 

For an algebra $A$, consider a right $\times_A$-Hopf algebra $\cC$, with
source map $s$, target map $t$, coproduct $\DC$, counit $\eC$ and
bijective canonical map
$$
\theta:\cC\stac {A^{op}} \cC\to \cC\stac A \cC,\qquad c\stac {A^{op}}
c'\mapsto c\DC(c')=: c c^{\prime(1)}\stac A c^{\prime(2)},
$$
as in \eqref{eq:can_x_A}. Let $P$ be an $A^{op}$-subring in $\cC$, i.e. a
subalgebra such that $t(a)\in P$, for all $a\in A$. Assume that $\DC(p)\in
\cC\sstac A P$ for all $p\in P$ (e.g. $P$ is a left subcomodule of $\cC$).
Denote by $P^+$ the intersection of $P$ with the kernel of $\eC$. Then the
right ideal (hence $A$-$A$ sub-bimodule) $P^+ \cC$, generated by $P^+$, is also
a coideal. Indeed, in a right bialgebroid $\ker \eC$ is a right ideal, hence
for $p\in P^+$ and $c\in \cC$ we have  
$
\eC(pc)=\eC\big( s(\eC(p)) c)=0.
$
Furthermore, using  the Takeuchi axiom for a right
bialgebroid in the second equality, one computes 
\begin{eqnarray}\label{eq:D_pc}
\DC(pc)&=&p^{(1)} c^{(1)}\stac A \big(p^{(2)}-t(\eC(p^{(2)}))\big)
c^{(2)}  + p^{(1)} c^{(1)}\stac A t(\eC(p^{(2)}))c^{(2)}\nonumber\\
&=& p^{(1)} c^{(1)}\stac A \big(p^{(2)}-t(\eC(p^{(2)}))\big)
c^{(2)}  + p^{(1)} s(\eC(p^{(2)}))c^{(1)}\stac A c^{(2)}\nonumber\\
&=& p^{(1)} c^{(1)}\stac A \big(p^{(2)}-t(\eC(p^{(2)}))\big)
c^{(2)}  + p c^{(1)}\stac A c^{(2)}.
\end{eqnarray}
The second term in \eqref{eq:D_pc} is a (possibly zero) element of
$P^+\cC\sstac A \cC$. The map $c\mapsto c-t(\eC( c))$ splits the inclusion
$\mathrm{ker}\ \eC\subseteq \cC$, and its restriction splits the inclusion
$P^+\subseteq P$. Hence the first
term in \eqref{eq:D_pc} belongs to $\cC\sstac A P^+\cC$. Thus $P^+\cC$ is a
coideal. Denote by $\cQ:=\cC/P^+\cC$ the quotient coring and right
$\cC$-module. Let $\pi:\cC\to \cQ$ denote the canonical epimorphism, which is a
homomorphism of $A$-corings and of right $\cC$-modules (so in particular of
left $A$-modules). The map $\pi$ induces a $\cQ$-comodule structure on $\cC$,
with coaction $\varrho^\cC:= (\cC\sstac A \pi)\circ \DC$. Denote the algebra of
$\cQ$-coinvariants in $\cC$ by $B$. Note that, similarly to the computations in
\eqref{eq:D_pc}, for all $p\in P$ and $c\in \cC$,
\begin{eqnarray*}
\varrho^\cC(pc) &= &p^{(1)} c^{(1)}\stac A
\pi\big(\big(p^{(2)}-t(\eC(p^{(2)}))\big) c^{(2)} \big) + p^{(1)}
c^{(1)}\stac A \pi\big(t(\eC(p^{(2)})) c^{(2)} \big)\\
&=& pc^{(1)} \stac A \pi\big(c^{(2)}\big). 
\end{eqnarray*}
Hence we have a sequence of algebra inclusions $t(A)\subseteq P \subseteq
B\subseteq \cC$. We claim that $B\subseteq \cC$ is a $\cQ$-Galois extension,
that is, the canonical map  
\begin{equation}\label{eq:can_Qhs}
\mathrm{can}: \cC\stac B \cC\to \cC\stac A \cQ,\qquad c\stac B c'\mapsto
cc^{\prime(1)} \stac A \pi(c^{\prime (2)}),
\end{equation}
is bijective. Consider the composite map
\begin{equation}\label{eq:can_Qhs_inv}
\xymatrix{
\cC\stac A \cC \ar[r]^{\theta^{-1}}& \cC\stac {A^{op}} \cC \ar@{>>}[r]&
\cC\stac{B} \cC},
\end{equation}
where the rightmost arrow denotes the canonical epimorphism induced by
the algebra inclusion $A^{op}\cong t(A)\subseteq B$. We show that the map
\eqref{eq:can_Qhs_inv} factorises through $\cC\sstac A \cQ$. Introduce the index
notation $\theta^{-1}(1_\cC\sstac A c)=:c_{-}\sstac {A^{op}} c_{+}$, where
implicit summation is understood. Since $\theta$ is a right $\cC$-comodule
map, so is $\theta^{-1}$. That is, for all $c\in \cC$, 
$$c_{-}\sstac{A^{op}}
{c_{+}}^{(1)} \sstac A {c_{+}}^{(2)} = {c^{(1)}}_{-} \sstac{A^{op}}
{c^{(1)}}_{+} \sstac A c^{(2)}.
$$
 This implies that, for all $p\in P$,
$$
p_{-} \stac{A^{op}} p_{+} = p_{-}  \stac{A^{op}} {p_{+}}^{(2)}
t\big(\eC({p_{+}}^{(1)} )\big) = {p^{(1)}}_{-} \stac{A^{op}} p^{(2)} 
t\big(\eC({p^{(1)}}_{+} )\big) 
$$
is an element of $\cC\sstac {A^{op}} P$, so a (possibly zero) element of
$\cC\sstac {A^{op}} B$.  
Hence, for all $c,c'\in 
\cC$ and $p\in P^+$, the left $\cC$-module map \eqref{eq:can_Qhs_inv} takes
$c'\sstac A pc$ to
$$
c'(pc)_{-} \stac B (pc)_{+} = c' c_{-} p_{-} \stac B p_{+} c_{+} = c' c_{-}
p_{-}p_{+} \stac B c_{+} = c' c_{-} s\big(\eC(p)\big) \stac B c_{+} =
0. 
$$
In the first equality we used the multiplicativity of the map
$\theta^{-1}(1_\cC\sstac A \bullet)$, i.e.\ that for all $c,c'\in \cC$,
 $(cc')_{-} \sstac {A^{op}}
(cc')_{+} = c'_{-} c_{-} \sstac {A^{op}} c_{+} c'_{+}.
$
 In
the third equality we used that $c_{-} c_{+} = s\big(\eC(c)\big)$. 
The last equality follows since $p\in P^+\subseteq
\ker \eC$. Thus we proved the existence of a map
$$
\mathrm{can}^{-1}: \cC\stac A \cQ\to \cC\stac B \cC,\qquad c\stac A
\pi(c')\mapsto cc'_{-}\stac B c'_{+}.
$$
Since it is defined in terms of $\theta^{-1}$, it is straightforward to see
that it is the inverse of the canonical map \eqref{eq:can_Qhs}.

By Lemma~\ref{lemma.tor1}, corresponding to the $\cQ$-Galois extension
$B\subseteq \cC$, there is an $A$-$B$ pre-torsor structure on $\cC$, with
pre-torsor map 
$$
\cC\to \cC\stac A \cC\stac B \cC,\qquad c\mapsto c^{(1)}\stac A
{c^{(2)}}_{-}\stac B {c^{(2)}}_{+}.
$$
\end{example}

\begin{example}
Consider a right entwining structure $(T,\cC,\psi)$ over an algebra $A$, such
that the right regular $T$-module 
extends to an entwined module. That is, $T$ is a right $\cC$-comodule, with
coaction $\varrho^T:T\to T\sstac A \cC$, $t\mapsto t^{[0]}\sstac A t^{[1]}$
such  that, for all $t,t'\in T$, 
\begin{equation} \label{eq:T_entw_m}
(tt')^{[0]} \stac A (tt')^{[1]} = t^{[0]}\psi(t^{[1]}\stac A t').
\end{equation}
Assume that the range of the unit map $\alpha:A\to T$ of the $A$-ring $T$ lies
within the coinvariant subalgebra $B:=T^{co\cC}$. Furthermore, assume that
there exists a left $A$-module right $\cC$-comodule map $j:\cC\to T$, which
possesses a convolution inverse, i.e. an $A$-$A$ bimodule map
$\widetilde{j}:\cC\to T$, such 
that  
$$
j(c^{(1)})\widetilde{j}(c^{(2)}) =\alpha\big(\eC(c)\big)=
\widetilde{j}(c^{(1)})j(c^{(2)}), \qquad \textrm{for all } c\in \cC,
$$
where $\eC$ denotes the counit of $\cC$, and $\DC(c) =
c^{(1)}\otimes_Ac^{(2)}$. In this situation $T$ is said to
be a {\em $\cC$-cleft extension} of $B$, cf. \cite[Proposition
  6.4]{BohmVer:Mor&cleft}. A $\cC$-cleft extension is $\cC$-Galois by
\cite[Corollary 5.3]{BohmVer:Mor&cleft}. Hence by Lemma \ref{lemma.tor1} there
is a corresponding $A$-$B$ pre-torsor structure on $T$. The pre-torsor map is
$$
\tau:T\to T\stac A T \stac B T,\qquad t\mapsto t^{[0]}\stac A
\widetilde{j}(t^{[1]}) \stac B j(t^{[2]}).
$$
Note that $\tau$ is well defined since $j$ and $\widetilde{j}$ are $A$-$A$
bimodule maps and $B$ is an $A$-ring. It is obviously $B$-$A$
bilinear. Also,
$$
(T\stac A \mu_T)(\tau(t))=
t^{[0]}\stac A \widetilde{j}(t^{[1]})j(t^{[2]})=
t^{[0]}\stac A \alpha\big(\eC (t^{[1]})\big)=
t^{[0]} \alpha\big(\eC (t^{[1]})\big) \stac A 1_T=
t \stac A 1_T,
$$
hence axiom (b) in Definition \ref{def.pre.tor} holds. Axiom (a) is proven as
follows. Similarly to \cite[Lemma 4.7 1]{Abu:Mor.cor}, the identity
\begin{equation}\label{eq:Abuh}
\widetilde{j}(c) 1_T^{[0]} \stac A  1_T^{[1]} =\psi(c^{(1)}\stac A
\widetilde{j}(c^{(2)}))
\end{equation} 
can be proven, for all $c\in \cC$. Hence, for all $t\in T$,
$$
\varrho\big( t^{[0]}\widetilde{j}(t^{[1]})\big) =
t^{[0]}\psi\big(t^{[1]}\stac A \widetilde{j}(t^{[2]})\big) =
t^{[0]}\widetilde{j}(t^{[1]})1_T^{[0]}\stac A 1_T^{[1]}.
$$
The first equality follows by \eqref{eq:T_entw_m} and the second one follows by
\eqref{eq:Abuh}. That is, $t^{[0]}\widetilde{j}(t^{[1]})\in B$, for all $t\in
T$. Hence
\begin{eqnarray*}
(\mu_T\stac B T)(\tau(t)) &=&
t^{[0]}\widetilde{j}(t^{[1]})\stac B j(t^{[2]})=
1_T \stac B t^{[0]}\widetilde{j}(t^{[1]})j(t^{[2]})\\
&=&
1_T \stac B t^{[0]}\alpha\big(\eC (t^{[1]})\big)=
1_T \stac B t.
\end{eqnarray*}
Finally, axiom (c) follows by the right $\cC$-colinearity of $j$ as
\begin{eqnarray*}
(T\stac A T\stac B \tau)(\tau(t))&=&
t^{[0]}\stac A \widetilde{j}(t^{[1]}) \stac B j(t^{[2]})^{[0]}\stac A 
\widetilde{j}(j(t^{[2]})^{[1]})\stac B j(j(t^{[2]})^{[2]}) \\
&=&t^{[0]}\stac A \widetilde{j}(t^{[1]}) \stac B j(t^{[2]})\stac A 
\widetilde{j}(t^{[3]})\stac B j(t^{[4]})=
(\tau \stac A T \stac B T)(\tau(t)).
\end{eqnarray*}
\end{example}

\section{Equivalences induced by pre-torsors}
\label{sec.equpretor}

Faithfully flat Hopf-Galois extensions play a central role in the description
of $k$-linear monoidal equivalences between comodule categories of flat Hopf
algebras 
over a commutative ring $k$. Extending results on commutative Hopf algebras in
\cite{SaaR:CatTan} and \cite{DelMil:TanCat}, in \cite{Ulb:Galfun} and
\cite{Ulb:Fibre}  Ulbrich established an equivalence of fibre functors
$U:{}^\cC\cM\to \cM_k$ with faithfully flat $\cC$-Galois extensions
$U(\cC)$ of $k$, for any $k$-flat Hopf algebra $\cC$. If $k$ is a field, then a
reconstruction theorem in \cite{Ulb:Rec} implies that for any fibre functor
$U$ there exists a (unique up to an isomorphism) $k$-Hopf algebra $\cD$, such
that 
$U$ factorises through a $k$-linear monoidal equivalence functor ${}^\cC\cM\to
{}^\cD \cM$ and the forgetful functor ${}^\cD\cM\to \cM_k$. Thus there is an
equivalence between $k$-linear monoidal equivalence functors ${}^\cC\cM\to
{}^\cD \cM$ and $\cC$-Galois extensions of $k$.

If $k$ is any commutative ring, then no reconstruction theorem for arbitrary
fibre functors is available. However, as Schauenburg pointed out in
\cite{Sch:big} and \cite{Sch:Hop}, for any faithfully flat $\cC$-Galois
extension $T$ of $k$ there exists a (unique up to an isomorphism) $k$-Hopf
algebra 
$\cD$, such that $T$ is a $\cD$-$\cC$ bi-Galois extension of $k$. This
observation was used in \cite{Sch:big} to prove that the fibre functor
${}^\cC\cM\to \cM_k$, induced by $T$, factorises through a $k$-linear
monoidal equivalence functor ${}^\cC\cM\to {}^\cD \cM$ and the forgetful
functor ${}^\cD\cM\to \cM_k$. Furthermore, every $k$-linear monoidal
equivalence functor ${}^\cC\cM\to {}^\cD \cM$ was shown to be induced by a
faithfully flat $\cD$-$\cC$ bi-Galois extension of $k$.

The first steps to extend the above theory to non-commutative base algebras
(replacing the commutative ring $k$ above), were made by Schauenburg in
\cite{Sch:Bianc}. 
Faithfully flat ${\mathcal C}$-Galois extensions of an arbitrary algebra $B$,
for a $k$-flat Hopf algebra $\cC$, were proven to induce $k$-linear monoidal  
equivalences between the categories of comodules of the $k$-Hopf algebra $\cC$
and a {\em $B$-bialgebroid} $\cD$.

Our aim in the next two sections is to develop a more symmetric study of
functors induced by faithfully flat (bi-)Galois extensions. That is, to replace
both Hopf algebras $\cC$ and $\cD$ above with bialgebroids over arbitrary, 
different base algebras $A$ and $B$. The problem is divided into two
parts. While the bialgebroid structures of $\cC$ and $\cD$ are essential for
having a monoidal structure on the categories of their comodules,
one can study more general functors
between categories of comodules of corings, induced 
by faithfully flat $A$-$B$ pre-torsors. By Theorem \ref{thm.main.1}, a
faithfully flat $A$-$B$ pre-torsor $T$ is a bi-Galois object, for a $B$-coring
$\cD$ and an 
$A$-coring $\cC$. In this section it is shown that (under some additional
assumptions (iii) in Remark \ref{rem:ff} below) the cotensor
product $T\Box_\cC\, \bullet$ defines a functor ${}^\cC\cM\to {}^\cD\cM$.
As a main result we prove that this functor
is an equivalence of categories. In forthcoming Section \ref{sec.tors} our
study will be 
specialised to functors induced by $A$-$B$ torsors, whose monoidal properties
will be studied.

\smallskip

In a similar way as in \cite[Theorem 2.7]{BrzHaj:coa}, to a right
$\cC$-Galois extension $B\subseteq T$ one associates a right entwining
structure over $A$, consisting of the $A$-ring $T$, the $A$-coring $\cC$ and
the entwining map
$$
\cC\stac A T \to T\stac A \cC,\qquad c\stac A t\mapsto 
\can_\cC\big(\can_\cC^{-1}(1_T\stac A c)\ t\big),
$$
such that the right regular $T$-module extends to a right entwined
module.
In a symmetric way 
to a left $\cD$-Galois extension $A\subseteq T$ there
corresponds a left entwining structure over $B$, consisting of the $B$-ring
$T$, the $B$-coring $\cD$ and the entwining map
$$
T\stac B \cD \to \cD\stac B T,\qquad t\stac B d\mapsto 
{}_\cD\can\big(t\  {}_\cD\can^{-1}(d\stac B 1_T)\big),
$$
such that the left regular $T$-module extends to a left entwined
module.

By Theorem \ref{thm.main.1}, to a faithfully flat $A$-$B$ pre-torsor $T$ one
can 
associate a right entwining structure $(T,\cC,\psi_\cC)$ over $A$ and a left
entwining structure $(T,\cD,\psi_\cD)$ over $B$ in the ways described
above. The corings $\cC$ and $\cD$ were determined in Theorem~\ref{thm.main.1}
(ii) and (iii), respectively. The entwining maps have the following explicit
forms. 
\begin{eqnarray}
\label{eq:psiC}
\psi_\cC:\cC\stac A T&\to T\stac A \cC,\qquad
\sum_i t_i\stac B u_i\stac A v &\mapsto 
\sum_i t_i\tau(u_i v),\\
\label{eq:psiD}
\psi_\cD: T\stac B \cD &\to \cD\stac B T,\quad
t\stac B \sum_j u_j \stac A v_j &\mapsto 
\sum_j \tau(tu_j)v_j.
\end{eqnarray}
By an easy extension of
  \cite[32.8~(2)]{BrzWis:cor} to a non-commutative base algebra $A$,
$\cC\sstac A T$ is a right entwined module for the right
  entwining structure $(T,\cC,\psi_\cC)$ over the algebra $A$, with right
  $T$-action $\cC\sstac A \mu_T$ and right 
  $\cC$-coaction $(\cC\sstac A \psi_\cC)\circ (\Delta_\cC\sstac A T)$.
Symmetrically, $T\sstac B \cD$ is a left entwined module for the left
  entwining 
  structure $(T,\cD,\psi_\cD)$ over the algebra $B$, 
with left $T$-action $\mu_T\sstac B \cD$ and left
  $\cD$-coaction $(\psi_\cD\sstac B \cD)\circ (T\sstac B \Delta_\cD)$.

\begin{proposition}
\label{prop:barT}
Let $T$ be a faithfully flat $A$-$B$ pre-torsor and $\cC$ and $\cD$ the
  associated $A$- and $B$-corings, respectively. The following $A$-$B$
  sub-bimodules of $T\sstac B T\sstac A T$ coincide:
\begin{rlist}
\item The coinvariants of the left entwined module $T\sstac B \cD$ for
  the left entwining structure $(T,\cD,\psi_\cD)$ over the algebra $B$;
\item The coinvariants of the right entwined module $\cC\sstac A T$ for
  the right entwining structure $(T,\cC,\psi_\cC)$ over the algebra $A$;
\item The intersection of $\cC\sstac A T$ and $T\sstac B \cD$.
\end{rlist}
The $A$-$B$ bimodule in $(i)$-$(iii)$ will be denoted by ${\bar T}$ in what
follows. 
\end{proposition}
\begin{proof}
We prove the equality of the bimodules (i) and (iii). The equality of
(ii) and (iii) follows by a symmetrical reasoning. By the form
\eqref{eq:psiD} of the entwining map $\psi_\cD$, the left $\cD$-coaction 
$T\sstac B \cD\to \cD\sstac B T\sstac B \cD$ comes out explicitly as
$$
\sum_j t_j\stac B u_j\stac A v_j \mapsto \sum_j 
(t_j u_j\sut 1)\sut 1\stac A(t_j u_j\sut 1)\sut 2\stac B (t_j u_j\sut 1)\sut 3
u_j\sut 2\stac B u_j\sut 3 \stac A v_j.
$$
If the element $\sum_j t_j\sstac B u_j\sstac A v_j\in T\sstac B \cD$ belongs
to $\cC\sstac A T$, i.e. the kernel of the map $\omega\sstac A T$ in
Lemma \ref{lemma.pure.eq}, then 
\begin{eqnarray}
\label{eq:TDco}
&&\sum_j (t_j u_j\sut 1)\sut 1\stac A(t_j u_j\sut 1)\sut 2\stac B (t_j u_j\sut
1)\sut 3 u_j\sut 2\stac B u_j\sut 3 \stac A v_j\\
&&\hspace{2.5in}
=\sum_j 1_T\sut 1\stac A 1_T \sut 2\stac B 1_T\sut 3 t_j\stac B u_j\stac A
v_j, \nonumber
\end{eqnarray}
hence $\sum_j t_j\sstac B u_j\sstac A v_j$ is a coinvariant in $T\sstac B
\cD$. Conversely,
assume that $\sum_j t_j\sstac B u_j\sstac A v_j$ is a coinvariant in $T\sstac B
\cD$, i.e. satisfies \eqref{eq:TDco}. Applying $T\sstac A\mu_T\sstac B T\sstac
A T$ to both sides of \eqref{eq:TDco} and using Definition \ref{def.pre.tor}
(b), we conclude that $\sum_j t_j\sstac B u_j\sstac A v_j$ lies in the kernel
of $\omega\sstac A T$, that is $\cC\sstac A T$. 
\end{proof}

\begin{lemma}
\label{lem:barTcom}
Let $T$ be a faithfully flat $A$-$B$ pre-torsor with structure map $\tau$, and
  $\cC$ and $\cD$ the associated $A$- and $B$-corings, respectively. The
  $A$-$B$ bimodule ${\bar T}$ in Proposition \ref{prop:barT} is a $\cC$-$\cD$
  bicomodule. Both the left $\cC$-coaction and the right $\cD$-coaction are
  given by the restriction of $T\sstac B \tau\sstac A T$.
\end{lemma}
\begin{proof}
Since $(\tau\sstac A T)(\cD)\subseteq \cD\sstac B \cD$ and ${\bar T}\subseteq
T\sstac B \cD$, it follows that $(T\sstac B\tau\sstac A T)({\bar T})\subseteq
T\sstac B \cD\sstac B \cD$. Furthermore, by Definition \ref{def.pre.tor} (c),
for  
$\sum_j t_j\sstac B u_j\sstac A v_j\in {\bar T}\subseteq \cC\sstac A T$,
\begin{eqnarray*}
\label{eq:TDDco}
&&\hspace{-.5cm}\sum_j (t_j{u_j\sut 1}\sut 1)\sut 1 \stac A  (t_j{u_j\sut
  1}\sut 1)\sut 2\stac B   
(t_j{u_j\sut 1}\sut 1)\sut 3{u_j\sut 1}\sut 2 \sstac B {u_j\sut 1}\sut 3 \stac
A {u_j\sut 2} \stac B u_j\sut 3\stac A v_j\\
&&=\sum_j (t_j{u_j\sut 1})\sut 1 \stac A  (t_j{u_j\sut 1})\sut 2\stac B  
(t_j{u_j\sut 1})\sut 3{u_j}\sut 2 \sstac B {u_j\sut 3}\sut 1 \stac A {u_j\sut
  3}\sut 2 \stac B {u_j\sut 3}\sut 3\stac A v_j\\
&&=\sum_j 1_T\sut 1\stac A 1_T\sut 2\stac B 1_T\sut 3 t_j\stac B u_j\sut
1\stac A u_j\sut 2\stac B u_j\sut 3 \stac A v_j.
\end{eqnarray*}
That is, $\sum_j t_j\sstac B \tau(u_j)\sstac A v_j$ is 
a coinvariant in the left entwined module $T\sstac B \cD\sstac B \cD$,
for the left entwining structure $(T,\cD,\psi_\cD)$. By a left handed version
of Lemma
\ref{lem:purecoC}, ${}^{co\cD} (T\sstac B \cD\sstac B \cD)={}^{co\cD} (T\sstac
B \cD)\sstac B \cD$, hence  $(T\sstac B\tau\sstac A T)({\bar T})\subseteq
{\bar T}\sstac B\cD$. The map $T\sstac B \tau\sstac A T$ is obviously right
$B$-linear. Its counitality and coassociativity follow by Definition
\ref{def.pre.tor} (b) and (c), respectively. This proves that ${\bar T}$ is a
right $\cD$-comodule, with coaction given by the restriction of $T\sstac B
\tau\sstac A T$. By a symmetrical reasoning ${\bar T}$ is also a left
$\cC$-comodule with coaction given by the restriction of $T\sstac B
\tau\sstac A T$. The commutativity of the left $\cC$-coaction and the right
$\cD$-coaction in ${\bar T}$ follows by Definition \ref{def.pre.tor} (c).
\end{proof}
\begin{corollary}
\label{cor:barTeq}
Let $T$ be a faithfully flat $A$-$B$ pre-torsor and $\cC$ and $\cD$ the
  associated $A$- and $B$-corings, respectively. Let ${\bar T}$ be the 
$\cC$-$\cD$ bicomodule in Lemma \ref{lem:barTcom}.
\begin{zlist}
\item $T\sstac B \cD$ and $T\sstac A {\bar T}$ are isomorphic as
  left $T$-modules left $\cD$-comodules and right $\cD$-comodules.
  
\item $\cC\sstac A T$ and ${\bar T}\sstac B T$ are isomorphic as
  right $T$-modules right $\cC$-comodules and left $\cC$-comodules.  
\end{zlist}
\end{corollary}
\begin{proof}
By the Galois Coring Structure Theorem \cite[28.19 (2) $(a)\Rightarrow
  (c)$]{BrzWis:cor}, the `coinvariants' functor ${}_T^\cD \cM(\psi_\cD)\to
  {}_A \cM$ is 
  an equivalence, with  inverse $T \sstac A \bullet$. Hence the counit
  of the adjunction is an isomorphism of left $T$-modules and left
  $\cD$-comodules. In particular, by characterisation of ${\bar T}$ in
  Proposition \ref{prop:barT} (i),
\begin{equation}
\label{eq:Dcou}
T\stac A {\bar T}\to T\stac B \cD,\qquad t'\stac A \sum_j t_j\stac B u_j\stac
A v_j\mapsto  \sum_j t' t_j\stac B u_j\stac A v_j,
\end{equation}
is an isomorphism of left $T$-modules and left $\cD$-comodules. Its inverse
has the form
\begin{equation}
\label{eq:Dcouinv}
T\stac B \cD\to T\stac A {\bar T},\qquad t\stac B \sum_j u_j\stac A v_j
\mapsto \sum_j t u_j\sut 1 \stac A u_j\sut 2\stac B u_j\sut 3 \stac A v_j.
\end{equation}
Since
\eqref{eq:Dcou} is obviously colinear with respect to the right
$\cD$-comodule structures of $T\sstac A {\bar T}$ and $T\sstac B \cD$, defined
via their second factors, this completes the proof of claim (1). Part (2) is
proven symmetrically. 
\end{proof}
\begin{theorem}
\label{thm:TcotT}
Let $T$ be a faithfully flat $A$-$B$ pre-torsor and $\cC$ and $\cD$ the
  associated $A$- and $B$-corings, respectively. Let ${\bar T}$ be the
  $\cC$-$\cD$ bicomodule in Lemma~\ref{lem:barTcom}. Then 
\begin{zlist}
\item $T \Box_\cC {\bar T}$ is a $\cD$-$\cD$ bicomodule via the 
$\cD$-coactions of $T$ and $\bar{T}$, and $T \Box_\cC {\bar T}\cong \cD$ as
$\cD$-$\cD$ bicomodules.  
\item ${\bar T}\Box_\cD T$ is a $\cC$-$\cC$ bicomodule via 
the $\cC$-coactions of $\bar{T}$ and $T$, and ${\bar T}\Box_\cD T\cong \cC$ 
as $\cC$-$\cC$ bicomodules. 
\end{zlist}
\end{theorem}
\begin{proof}
Consider a map
\begin{equation}
\label{eq:TbTtoD}
\varpi:T\cstac\cC{\bar T}\to \cD,\qquad
\sum_j t_j\stac A u_j\stac B v_j\stac A w_j \mapsto
\sum_j t_j u_j v_j \stac A w_j =\sum_j t_j \stac A u_j v_j w_j,
\end{equation}
where the last equality follows by $\sum_j t_j\sstac A u_j\sstac B v_j\sstac A
w_j\in T\Box_\cC{\bar T}\subseteq T\sstac A
\cC\sstac A T$, hence $\sum_j t_j \sstac A u_j v_j \sstac A w_j\in T\sstac A A
\sstac A T$. By the defining equaliser property of the cotensor product 
$T\Box_\cC{\bar T}$, Definition \ref{def.pre.tor} (a), and since ${\bar
  T}\subseteq 
T\sstac B \cD$,  
\begin{eqnarray*}
\sum_j t_j\sut 1\stac A t_j\sut 2 \stac B t_j\sut 3 u_j v_j w_j 
&= &\sum_j t_j\stac A u_j \stac B v_j\sut 1 v_j\sut 2 v_j\sut 3 w_j\\
&=&
\sum_j t_j \stac A u_j \stac B v_j w_j=
\sum_j t_j \stac A u_j v_j w_j \stac B 1_T.
\end{eqnarray*}
Hence the range of the map \eqref{eq:TbTtoD} is in $\cD$, indeed. Since
$\tau(T)$ 
lies within $T\sstac A \cC$, it follows that $(\tau\sstac A T)(\cD)$ lies
within $T\sstac A \cC\sstac A T$. On the other hand, $(\tau\sstac A T)(\cD)$
lies within $\cD\sstac B \cD\subseteq T\sstac A T\sstac B \cD$, hence
characterisation of ${\bar T}$ in Proposition \ref{prop:barT}~(iii), together
with the flatness of $T$ as a right $A$-module, implies that 
$(\tau\sstac A T)(\cD)$ is a subset of $T\sstac A
{\bar T}$. Furthermore, by Definition~\ref{def.pre.tor}~(c), the restriction
of $\tau\sstac A T$ defines a map $\cD\to T\Box_\cC {\bar T}$. We claim that
this map is the inverse of $\varpi$ in \eqref{eq:TbTtoD}. By 
Definition~\ref{def.pre.tor}~(b) (or (a)), the restriction of $\varpi \circ
(\tau\sstac A T)$ to 
$\cD$ is equal to the identity map in $\cD$. On the other hand, for $x = \sum_j
t_j\sstac A u_j \sstac B v_j \sstac A w_j\in T\Box_\cC {\bar T}$, 
$$
(\tau\stac A T)\big(\varpi(x)\big)
\!=\!\sum_j t_j\sut 1 \stac A t_j\sut 2 \stac B t_j \sut 3 u_j v_j \stac A w_j
\!=\! \sum_j t_j \stac A u_j \stac B v_j \sut 1 v_j \sut 2 v_j \sut 3 \stac A
w_j  
\!=\! x.
$$
Here the second equality follows by the definition of the cotensor product
$T\Box_\cC {\bar T}$ and the third one follows by Definition
\ref{def.pre.tor} (b). This proves that the map \eqref{eq:TbTtoD} is
bijective. Next take any 
$x= \sum_j t_j\sstac A u_j \sstac B v_j \sstac A w_j\in T\Box_\cC {\bar T}$
and compute
\begin{eqnarray*}
(\cD\stac B \varpi^{-1})\big(\Delta_\cD(\varpi(x))\big)
&=&\sum_j t_j\sut 1\stac A t_j\sut 2\stac B {t_j\sut 3}\sut 1\stac A
{t_j\sut 3}\sut 2\stac B {t_j\sut 3}\sut 3 u_j v_j \stac A w_j\\
&=& \sum_j {t_j\sut 1}\sut 1\stac A {t_j\sut 1}\sut 2\stac B  {t_j\sut 1}\sut
3\stac A {t_j}\sut 2\stac B {t_j}\sut 3 u_j v_j \stac A w_j\\
&=& \sum_j {t_j}\sut 1\stac A {t_j}\sut 2\stac B  {t_j}\sut 3\stac A
u_j\stac B {v_j}\sut 1 v_j\sut 2 v_j\sut 3 \stac A w_j\\
&=& \sum_j {t_j}\sut 1\stac A {t_j}\sut 2\stac B  {t_j}\sut 3\stac A
u_j\stac B {v_j}\stac A w_j
=(\tau\stac A {\bar T})(x).
\end{eqnarray*}
The second equality follows by Definition \ref{def.pre.tor} (c), the third one
follows by the equaliser property of the cotensor product $T\Box_\cC {\bar T}$
and the fourth one does by Definition~\ref{def.pre.tor}~(b). The above
computation confirms that $T\Box_\cC{\bar T}$ is a left $\cD$-comodule
with the coaction given by the restriction of $\tau\ot_A\bar{T}$ and that
$\varpi$ is left $\cD$-colinear. The fact that $T\Box_\cC{\bar T}$ is a 
right $\cD$-comodule and the right
$\cD$-colinearity of $\varpi$  are proven by a similar computation, for 
$x= \sum_j t_j\sstac A u_j \sstac B v_j \sstac A w_j\in T\Box_\cC {\bar T}$.
\begin{eqnarray*}
(\varpi^{-1}\stac B \cD)\big(\Delta_\cD(\varpi(x))\big)
&=& \sum_j {t_j\sut 1}\sut 1\stac A {t_j\sut 1}\sut 2\stac B  {t_j\sut 1}\sut
3\stac A {t_j}\sut 2\stac B {t_j}\sut 3 u_j v_j \stac A w_j\\
&=&\sum_j t_j\sut 1\stac A t_j\sut 2\stac B {t_j\sut 3}\sut 1\stac A
{t_j\sut 3}\sut 2\stac B {t_j\sut 3}\sut 3 u_j v_j \stac A w_j\\
&=& \sum_j t_j\stac A u_j \stac B {v_j\sut 1}\sut 1 \stac A {v_j\sut 1}\sut 2
\stac B {v_j\sut 1}\sut 3 v_j\sut 2 v_j\sut 3\stac A w_j\\
&=& \sum_j t_j\stac A u_j \stac B {v_j}\sut 1 \stac A {v_j}\sut 2 \stac B
{v_j}\sut 3 \stac A w_j\\
&=&(T\stac A(T\stac B \tau\stac A T))(x).
\end{eqnarray*} 
Here again, the second equality follows by Definition \ref{def.pre.tor} (c),
the third one follows by the equaliser property of the cotensor product
$T\Box_\cC {\bar T}$ and the fourth one does by Definition \ref{def.pre.tor}
(b) and the right $A$-linearity of $\tau$. This completes the proof of part
(1). Part (2) follows by symmetrical arguments.
\end{proof}
\begin{lemma} 
\label{lem:flat}
  Let $T$ be a faithfully flat $A$-$B$ pre-torsor. 
\begin{zlist}
\item If $T$ is a flat right $B$-module, then the $A$-coring $\cC$, associated
  to $T$ in Theorem~\ref{thm.main.1}~(ii), is a flat right $A$-module.
\item If $T$ is a flat left $A$-module, then the $B$-coring $\cD$, associated
  to $T$ in Theorem~\ref{thm.main.1}~(iii), is a flat left $B$-module. 
\end{zlist}
\end{lemma}
\begin{proof}
    Recall that for any ring extension $A\subseteq T$, the corresponding
    restriction of scalars (forgetful) functor ${}_T \cM \to {}_A \cM$ is 
    faithful, hence it reflects monomorphisms. On the other hand, it
    possesses a left adjoint (the induction functor $T\sstac A\, \bullet$),
    hence it preserves monomorphisms as well.

    $T$ is a faithfully flat right $A$-module by assumption. Hence $\cC$ is a
    flat right $A$-module, i.e. the functor
    $\cC\sstac A \, \bullet:{}_A \cM\to \cM_k$ preserves 
    monomorphisms, if and only if the functor $T\sstac A \cC\sstac A \, \bullet
    :{}_A \cM \to \cM_k$ preserves monomorphisms. By the right
    $\cC$-Galois property of the 
    extension $B\subseteq T$ (cf.\ Theorem~\ref{thm.main.1}), the right
    $A$-modules $T\sstac A \cC$ and 
    $T\sstac B T$ are isomorphic. Hence the functors $T\sstac A \cC\sstac
    A \, \bullet :{}_A \cM \to \cM_k$ and $T\sstac B T\sstac A \, \bullet
    :{}_A \cM \to \cM_k$ are naturally isomorphic. The flatness of
    $\cC$ as a right $A$-module follows by the assumption that both functors
    $T\sstac A \,\bullet :{}_A \cM \to {}_B \cM$ and $T\sstac B \,
    \bullet:{}_B \cM \to \cM_k$ preserve monomorphisms, hence so
    does their composite. This completes the proof of claim (1).
    Assertion (2) follows by a symmetrical reasoning.
\end{proof}

In light of Corollary \ref{cor:barTeq}, analogous considerations to those used
to prove Lemma \ref{lem:flat}, lead to the following 
\begin{lemma}\label{lem:flat2}
Let $T$ be a faithfully flat $A$-$B$ pre-torsor. Assume that $T$ is 
  faithfully flat also as a right $B$-module. Then the associated $B$-coring
  $\cD$ is a flat right $B$-module if and only if the right $B$-module ${\bar
  T}$ in Proposition \ref{prop:barT} is flat.
\end{lemma}
\begin{remark}\label{rem:ff}
Consider an $A$-$B$ pre-torsor $T$ with structure map $\tau$. If both $A$ and
$B$ coincide with the ground ring $k$, then the following properties of $T$ are
equivalent by Lemmas \ref{lem:flat} and \ref{lem:flat2}. 
\begin{rlist}
\item $T$ is a faithfully flat pre-torsor, i.e.\ a faithfully flat left
  $B$-module and right $A$-module.
\item $T$ is a faithfully flat left and right $B$-module and right $A$-module.
\item $T$ is a faithfully flat left and right $B$-module and right $A$-module, 
and the left $B$-module map $\tau\sstac A M -T\sstac A {}^M\varrho$ is
$\cD\sstac B \cD$-pure, for any left $\cC$-comodule $M$, with coaction ${}^M
\varrho$.
\item $T$ is a faithfully flat left and right $B$-module and right $A$-module, 
and the associated $B$-coring in Theorem \ref{thm.main.1} is a flat right
$B$-module.
\item $T$ is a faithfully flat left and right $B$-module and right $A$-module, 
and the associated $A$-$B$ bimodule ${\bar T}$ in Proposition \ref{prop:barT} 
is a flat right $B$-module.
\end{rlist}
For arbitrary base algebras $A$ and $B$, properties (i)-(v) seem no longer 
equivalent. Only implications (i)$\Leftarrow$ (ii) $\Leftarrow$ (iii)
$\Leftarrow$ (iv)$\Leftrightarrow$ (v) are easily proven. Indeed,
(iv)$\Leftrightarrow$(v) is proven in Lemma \ref{lem:flat2}. If $\cD$ is a
flat right $B$-module then so is $\cD\sstac B \cD$. Hence any left $B$-module
map is $\cD\sstac B \cD$-pure. This proves (iv)$\Rightarrow$(iii). The
remaining implications (iii)$\Rightarrow$(ii)$\Rightarrow$(i) are obvious.
Some theorems can be
proven by assuming a weaker one of the above properties, in other cases a
stronger one is needed.
\end{remark}
\begin{corollary}\label{cor:prtoreq}
Let $T$ be an $A$-$B$ pre-torsor obeying properties (iii) in 
Remark~\ref{rem:ff}.  
Then the categories of left comodules of the associated $A$- and $B$-corings
  $\cC$ and $\cD$ are
  equivalent. The inverse equivalences between them are given by the cotensor
  products $T\Box_\cC\, \bullet: {}^\cC\cM\to {}^\cD \cM$ and ${\bar
  T}\Box_\cD\, \bullet :{}^\cD\cM\to {}^\cC \cM$, respectively, where ${\bar
  T}$ is the $\cC$-$\cD$ bicomodule in Lemma~\ref{lem:barTcom}.
\end{corollary}
\begin{proof}
By the $\cD\sstac B \cD$-purity of the equalisers defining  $\cC$-cotensor
products with $T$, there is a functor
$T\Box_\cC\, \bullet:{}^\cC\cM\to {}^\cD\cM$. By Lemma~\ref{lem:flat}~(1),
$\cC$ is a flat 
right $A$-module, so the $\cC$-$\cD$ bicomodule ${\bar T}$ in 
Lemma~\ref{lem:barTcom} induces another functor ${\bar T}\Box_\cD\,
\bullet:{}^\cD\cM\to {}^\cC\cM$. By a reasoning similar to \cite[Section
2]{Sch:big}, the following sequence of isomorphisms holds, for any left
$\cD$-comodule $M$,
\begin{eqnarray*}
T\stac B [(T\cstac \cC {\bar T})\cstac \cD M]
&\cong& [T\stac B (T\cstac \cC {\bar T})]\cstac \cD M
\cong [(T\stac B T)\cstac \cC {\bar T}]\cstac \cD M
\cong [(T\stac A \cC)\cstac \cC {\bar T}]\cstac \cD M\\
&\cong&(T\stac A {\bar T})\cstac \cD M
\cong T\stac A ({\bar T}\cstac \cD M)
\cong (T\stac A\cC)\cstac \cC ({\bar T}\cstac \cD M)\\
&\cong& (T\stac B T)\cstac \cC ({\bar T}\cstac \cD M)
\cong T\stac B [T\cstac \cC ({\bar T}\cstac \cD M)].
\end{eqnarray*}
The first two and the last isomorphisms follow by the flatness of $T$ as a
right $B$-module. The fifth one follows by the flatness of $T$ as a right
$A$-module. The third and the penultimate equivalences follow by the right
$\cC$-Galois property of the extension $B\subseteq T$ (cf.\ 
Theorem~\ref{thm.main.1}). 
Since $T$ is a faithfully flat right $B$-module by assumption, this implies
that $(T\Box_ \cC {\bar T})\Box_\cD M\cong T\Box_ \cC ({\bar T}\Box_ \cD M)$,
for any left 
$\cD$-comodule $M$. Together with Theorem \ref{thm:TcotT} (1), this
proves that the composite of the functors ${\bar T} \Box_\cD\,\bullet $ and
$T\Box_\cC \,\bullet$ is naturally isomorphic to the identity functor ${}^\cD
\cM$. By Corollary \ref{cor:barTeq} (1), an analogous argument proves that
$T\sstac 
A [({\bar T}\Box_\cD T)\Box_\cC N]\cong T\sstac A [{\bar T}\Box_\cD (T\Box_\cC
N)]$, for any left $\cC$-comodule $N$. Since $T$ is a faithfully flat right
$A$-module, we conclude that 
$({\bar T}\Box_\cD T)\Box_\cC N\cong {\bar T}\Box_\cD (T\Box_\cC N)$. 
Hence it follows by Theorem \ref{thm:TcotT} (2) that the composite of the
functors $T\Box_\cC \,\bullet$ and ${\bar T} \Box_\cD\,\bullet $ is naturally
isomorphic to the identity functor on ${}^\cC\! \cM$.
\end{proof}

In the 
rest of the section we study the particular case when the entwining maps
$\psi_\cC$ in \eqref{eq:psiC} and $\psi_\cD$ in \eqref{eq:psiD} are bijective. 
By the standard entwining structure arguments,
$(T,\cC,\psi_\cC^{-1})$ is a left entwining structure over $A$, and $T$ is a
left entwined module with the left regular $T$-action and the left
$\cC$-coaction 
\begin{equation}\label{eq:lcoac}
{}^T\!\varrho:T\to \cC\stac A T,\qquad
t\mapsto \psi_\cC^{-1}\big(t\tau(1_T)\big).
\end{equation}
The coinvariants of $T$ as a left $\cC$-comodule coincide with the
coinvariants $B$ of $T$ as a right $\cC$-comodule. The left $\cC$-Galois
property of the algebra extension $B\subseteq T$ is equivalent to its right
$\cC$-Galois property. The canonical maps are related by the entwining map,
i.e.
\begin{equation}\label{eq:cans}
\can_\cC=\psi_\cC\circ {}_\cC \can.
\end{equation}
By these considerations, if for a faithfully flat
$A$-$B$ pre-torsor $T$, the entwining maps \eqref{eq:psiC} and \eqref{eq:psiD}
are bijective, then $B\subseteq T$ is a left $\cC$-Galois extension and
$A\subseteq T$ is a right $\cD$-Galois extension. In fact one can prove more.
\begin{theorem}\label{thm:optor}
Let $T$ be a faithfully flat $A$-$B$ pre-torsor and $\cC$ and $\cD$ the
associated corings in Theorem \ref{thm.main.1}. Assume that the entwining
maps \eqref{eq:psiC} and \eqref{eq:psiD} are bijective. Then the $\cC$-$\cD$
bicomodule ${\bar T}$ in Lemma \ref{lem:barTcom} is isomorphic to $T$. 
Therefore $T$ is a $\cC$-$\cD$ bi-Galois object, hence in particular 
a $B$-$A$ pre-torsor. 
\end{theorem}
\begin{proof}
For any element $t\sstac B u\sstac A v$ of $T\sstac B T \sstac A T$,
\begin{eqnarray*}
(T\stac A \psi_\cC)\big((\can_\cC \stac A T)(t\stac B u\stac A v)\big)
&=& t u\sut 1\stac A u\sut 2\tau(u\sut 3 v)
= (T\stac A \can_\cC)(t \tau(u) v), \\
(\psi_\cD\stac B T)\big((T\stac B {}_\cD\can)(t\stac B u\stac A v)\big)
&=& \tau(t u\sut 1)
u\sut 2 \stac B u\sut 3 v
= ({}_\cD \can \stac B T) (t \tau(u) v).
\end{eqnarray*}
Therefore, by \eqref{eq:cans} and its analogue for the $\cD$-Galois extension
$A\subseteq T$, 
\begin{equation}
\label{eq:canid}
(T\stac A {}_\cC\can^{-1})\circ (\can_\cC\stac A T)=(\can_\cD ^{-1}\stac B
T)\circ (T\stac B {}_\cD \can).
\end{equation}
Taking the inverses of both sides of \eqref{eq:canid} and evaluating on an
element $1_T \sstac A t \sstac B 1_T$ of $T\sstac A T\sstac B T$, we conclude
that the left $\cC$-coaction and the right $\cD$-coaction in $T$ map 
an element $t\in T$ to the same element of $T\sstac B T \sstac A T$.
For these equal coactions we use the notation ${\bar \tau}: T\to T \stac B T
\stac A T$. 
By the counitality of the left $\cC$-coaction or the right $\cD$-coaction in
$T$, ${\bar \tau}$ is an injective $A$-$B$ bimodule map $T\to T\sstac B T
\sstac A T$. As its range lies both in $\cC\sstac A T$ and $T\sstac B \cD$, 
${\bar \tau}(T)\subseteq {\bar T}$. In order to prove the converse inclusion,
take $\sum_i u_i\sstac B v_i \sstac A t_i\in {\bar T}$ and compute 
\begin{eqnarray*}
{\bar \tau}(\sum_i u_i v_i t_i)
&=&\psi_\cC ^{-1}\big(\sum_i u_i v_i t_i \tau(1_T)\big)
=\psi_\cC ^{-1}\big(\sum_i u_i v_i\sut 1 v_i\sut 2 \tau(v_i\sut 3 t_i)\big)\\
&=&\psi_\cC ^{-1}\big(\sum_i u_i\tau(v_i t_i)\big)
=\sum_i u_i\stac B v_i \stac A t_i.
\end{eqnarray*}
The second equality follows by the property that $\sum_i u_i\sstac B v_i
\sstac A t_i$ is an element of ${\bar T}\subseteq T\sstac B \cD$, the third
one does by Definition \ref{def.pre.tor} (a) and the left $B$-linearity of
$\tau$, and the last one follows by the form \eqref{eq:psiC} of
$\psi_\cC$. Thus ${\bar \tau}$ is a bijection $T\to {\bar T}$.
Its left $\cC$-colinearity and
right $\cD$-colinearity follow by the coassociativity of the
left $\cC$-coaction ${\bar \tau}$ and the right $\cD$-coaction ${\bar \tau}$
in ${\bar T}$, respectively. 
\end{proof}

Note that it follows by Theorem~\ref{thm:optor} 
that if for a faithfully flat $A$-$B$ pre-torsor $T$ the maps \eqref{eq:psiC}
and \eqref{eq:psiD} are bijective, then properties (ii), (iii), (iv) and (v) 
in Remark \ref{rem:ff} are all equivalent to each other.

The content of the following lemma was observed in \cite[Remark 2.4
(1)]{BohBrz:str}. 
\begin{lemma}
\label{lem:bijentw}
Let $T$ be a faithfully flat $A$-$B$ pre-torsor.
\begin{zlist}
\item If $T$ is a faithfully flat right $B$-module and the entwining map
  \eqref{eq:psiC} is bijective, then also the entwining map \eqref{eq:psiD} is
  bijective.
\item If $T$ is a faithfully flat left $A$-module and the entwining map
  \eqref{eq:psiD} is bijective, then also the entwining map \eqref{eq:psiC} is
  bijective.
\end{zlist}
\end{lemma}
\begin{proof}
    In the case when the entwining map \eqref{eq:psiC} is bijective, $T\sstac
    A T$ is a left entwined module for the left 
    entwining structure $(T,\cC,\psi_\cC^{-1})$ over $A$, via the left regular
    $T$-module structure of the first factor and the left $\cC$-coaction
    $(\psi_\cC^{-1}\sstac A T)\circ (T\sstac A {}^T\!\varrho)$, given in terms
    of the left $\cC$-coaction \eqref{eq:lcoac} in $T$. Its coinvariants
    are the elements $\sum_j u_j \sstac A v_j$ of $T\sstac A T$, for which
\begin{equation}
\label{eq:Dco1}
\sum_j (\psi_\cC^{-1}\stac A T)\big((T\stac A {}^T\!\varrho)(u_j\stac A
v_j)\big) =
\sum_j \psi_\cC^{-1}(\tau(1_T))u_j \stac A v_j. 
\end{equation}
Applying $\psi_\cC\sstac A T$ to both sides of \eqref{eq:Dco1}, we obtain the
equivalent condition
\begin{equation}
\label{eq:Dco2}
(T\stac A {}^T\! \varrho)(\sum_j u_j\stac A v_j)=
(\tau\stac A T)(\sum_j u_j\stac A v_j).
\end{equation}
By Theorem \ref{thm.main.1}, $B\subseteq T$ is a right $\cC$-Galois
extension. Hence one can apply the isomorphism $T\sstac A \can_\cC^{-1}\circ
\psi_\cC=T\sstac A T\sstac B \mu_T$ to \eqref{eq:Dco2} to obtain
$$
(T\stac A T\stac B \beta)(\sum_j u_j\stac A v_j)=
(T\stac A T \stac B \mu_T)\big((\tau\stac A T)(\sum_j u_j\stac A v_j)\big).
$$
This shows that the coinvariants of the left entwined module $T\sstac A
T$ coincide with the elements of $\cD$. Since the $T$-coring $T\sstac A \cC$,
associated to the right entwining structure $(T,\cC,\psi_\cC)$, is a Galois 
coring, so is the isomorphic coring $\cC\sstac A T$, associated 
to the left entwining structure $(T,\cC,\psi_\cC^{-1})$. By the assumption
that $T$ is a faithfully flat right $B$-module, the Galois Coring Structure
Theorem \cite[28.19 (2) (a) $\Rightarrow$ (c)]{BrzWis:cor} implies the
bijectivity of the map
\begin{equation}
\label{eq:canD}
T\stac B \cD \to T\stac A T,\qquad t\stac B \sum_j u_j \stac A v_j \mapsto
\sum_j tu_j \stac A v_j,
\end{equation} 
which is simply the counit of the adjunction of the `coinvariants' functor
${}^\cC_T 
\cM\to {}_B \cM$ and the induction functor $T\sstac B\, \bullet : {}_B \cM \to
{}^\cC_T \cM$, evaluated at the left entwined module $T\sstac A T$.
The map \eqref{eq:canD} is the inverse of the right $\cD$-canonical map
$\can_\cD$. Hence it is related via the $\cD$-analogue of \eqref{eq:cans} to 
the left $\cD$-canonical map ${}_\cD \can: T\sstac A T
\to \cD\sstac B T$, which is bijective by Theorem~\ref{thm.main.1}. Hence 
\eqref{eq:psiD} is bijective, as stated in (1).  
Part (2) is proven by a symmetric argument.
\end{proof}

\begin{corollary}
\label{cor:eq}
Let $T$ be an $A$-$B$ pre-torsor, obeying properties (ii) in Remark
\ref{rem:ff}. Assume that the entwining map \eqref{eq:psiC} is
bijective. Then $T$ possesses a $\cD$-$\cC$ bicomodule structure (cf.\ 
Theorem~\ref{thm.main.1}) and a $\cC$-$\cD$ bicomodule structure (cf.\
Theorem~\ref{thm:optor}). Furthermore, the functors  
$T\Box_\cD \,\bullet :{}^\cD \cM \to {}^\cC \cM$ and $T\Box_\cC \,
\bullet :{}^\cC \cM \to {}^\cD \cM$ are inverse equivalences.
\end{corollary}

\section{$A$-$B$ torsors as monoidal functors}
\label{sec.tors} 

In this section we focus our attention on faithfully flat $A$-$B$ torsors in
the sense of \cite[Definition 5.2.1]{Hobst:PhD}. Following \cite[Theorem
5.2.10]{Hobst:PhD}, in Theorem \ref{thm.main.Hobst} we establish a bijective
correspondence between faithfully flat $A$-$B$ torsors and bi-Galois objects
by bialgebroids (actually $\times_A$- and $\times_B$-Hopf algebras), $\cC$ and
$\cD$. 
A characteristic feature of bialgebroids is the monoidality of the category of
their comodules. In Theorems \ref{thm:lax} and \ref{thm:strmon} we show that an
$A$-$B$ torsor $T$, with properties in Remark~\ref{rem:ff}~(ii), induces a 
monoidal functor from ${}^\cC \cM$, the category of comodules of the associated
$\times_A$-Hopf algebra $\cC$, to the category of $B$-$B$ bimodules.
If $T$ obeys properties (iii) in Remark \ref{rem:ff}, then the induced functor
factorises through the category of $\cD$-comodules, monoidally.
What is more, by virtue of Corollary \ref{cor:prtoreq}, it results in a
monoidal equivalence between the comodule categories ${}^\cC \cM$ and ${}^\cD
\cM$ (cf.\ Theorem~\ref{thm:tor-moneq}).
In contrast to the case when $\cC$ and $\cD$ are (flat) Hopf algebras, it is
not known in general if all monoidal equivalence functors ${}^\cC\!\cM \to
{}^\cD \!\cM$ arise in this way.

We start by recalling the following \cite[Definition 5.2.1]{Hobst:PhD}.
\begin{definition}\label{def:ABtor}
An $A$-$B$ pre-torsor $T$ 
with unit maps $\alpha: A\to T$ and $\beta:B\to T$
 is called an {\em $A$-$B$ torsor} if 
$\alpha(A)$ and $\beta(B)$ are commuting subalgebras in $T$ and the 
structure map 
$\tau: T\to T\sstac A T\sstac B T$, $t\mapsto t\sut 1\sstac A t\sut 2\sstac B
t\sut 3$ obeys the following properties.
\begin{blist}
\item $\alpha(a)t\sut 1\stac A t\sut 2\stac B t\sut 3=t\sut 1 \stac A t\sut 2
  \alpha(a)\stac B   t\sut 3$; 
\item $t\sut 1 \stac A \beta(b) t\sut 2 \stac B t\sut 3=t\sut 1\stac A t\sut 2
  \stac B t\sut   3 
  \beta(b)$;
\item $\tau(t t')=t\sut 1 t^{\prime\langle 1 \rangle}\stac A 
t^{\prime\langle 2 \rangle} t\sut 2 \stac B 
t\sut 3 t^{\prime\langle 3 \rangle}$;
\item $\tau(1_T) = 1_T\stac{A} 1_T\stac{B} 1_T$,
\end{blist}
for all elements $t$ and $t'$ in $T$, $a$ in $A$ and $b$ in $B$.

An $A$-$B$ torsor is {\em faithfully flat} if it is faithfully flat as a right
$A$-module and left $B$-module.
\end{definition}
Note that axioms (a) and (b) are meaningful by the assumption that $\alpha(A)$
and $\beta(B)$ are commuting subalgebras in $T$. Axiom (c) makes sense in view
of axioms (a) and (b). In order to simplify
notation, we will not write out the unit maps $\alpha$ and $\beta$
explicitly in the sequel. 

The notion of an $A$-$B$ torsor is made interesting by its
relation to Galois extensions by bialgebroids.
It was observed in \cite[Theorem 5.2.10]{Hobst:PhD} that an
$A$-$B$ torsor $T$ determines a left, and a right Galois extension, by two
canonically associated bialgebroids, provided $T$ is 
faithfully flat as a left and right module for both base algebras $A$ and $B$
(cf.\ Remark~\ref{rem:rel_to_Hobst}). In contrast, in Theorem
\ref{thm.main.Hobst} below we assume faithful flatness of $T$ as a right
$A$-module and a left $B$-module only. We also prove the
converse of \cite[Theorem 5.2.10]{Hobst:PhD}, i.e.\ that a faithfully flat
(left or right) Galois extension by a bialgebroid determines a torsor. Putting 
these results together, we prove that the notions of a faithfully flat $A$-$B$
torsor, and that of a faithfully flat bi-Galois extension by bialgebroids, are
equivalent. For the convenience of the reader we include (in a sketchy form)
the complete proof, also of the parts which were obtained already in
\cite{Hobst:PhD}. Instead of following the original arguments there (operating
with the faithful flatness over its base algebra of a bialgebroid
associated to a torsor), we make use of Theorem~\ref{thm.main.1}. 
\begin{theorem} \label{thm.main.Hobst}
There is a bijective correspondence between the following sets of data:
\begin{rlist}
\item faithfully flat $A$-$B$ torsors $T$;
\item right $\times_A$-Hopf algebras $\cC$ and left faithfully flat right
  $\cC$-Galois  
  extensions $B\subseteq T$, such that $T$ is a right faithfully flat $A$-ring;
\item left $\times_B$-Hopf algebras $\cD$ and right faithfully flat left
  $\cD$-Galois 
  extensions $A\subseteq T$, such that $T$ is a left faithfully flat $B$-ring.
\end{rlist}
\end{theorem}
\begin{proof}
In light of Theorem \ref{thm.main.1}, we need to show that a faithfully flat
$A$-$B$ pre-torsor $T$ is a torsor if and only if the associated $A$-coring
$\cC$ is a right $\times_A$-Hopf algebra and $T$ is its right comodule
algebra. Equivalently, if and only if the associated $B$-coring
$\cD$ is a left $\times_B$-Hopf algebra and $T$ is its left comodule
algebra.

Assume first that $T$ is a faithfully flat $A$-$B$ torsor, with structure map
$\tau$. Use Definition \ref{def.pre.tor} (a), Definition~\ref{def:ABtor}~(b) 
and the definition of
$\cC$, as the kernel of the map $\omega$ in Lemma \ref{lemma.pure.eq}, to see
that, for $b\in B$ and $\sum_i u_i\sstac B v_i\in \cC$,
$$
\sum_i u_i \stac B v_i b=\sum_i u_i v_i \sut 1 v_i \sut 2 \stac B v_i \sut 3
b= \sum_i u_i v_i \sut 1 b v_i \sut 2 \stac B v_i \sut 3
= \sum_i bu_i\stac B v_i.
$$
This implies that, for any elements $\sum_i u_i\sstac B v_i$ and 
$\sum_j u'_j\sstac B v'_j$ in $\cC$, $\sum_{i,j} u'_j u_i\sstac B v_i v'_j$ is
a well defined element of $T\sstac B T$. Furthermore, it follows by
Definition~\ref{def:ABtor}~(c) that it belongs to $\ker \omega= \cC$. Since by 
Definition~\ref{def:ABtor}~(d) 
also $1_T\sstac B 1_T$ is an element of $\cC$, we conclude
that $\cC$ is an algebra, with multiplication and unit inherited from the
algebra $T^{op}\sstac k T$. The right $\cC$-coaction $\tau$ is multiplicative
by Definition~\ref{def:ABtor}~(c) and unital by
Definition~\ref{def:ABtor}~(d). Clearly, the maps
$$
A\to \cC, \quad a\mapsto 1_T \stac B a, \qquad \mbox{and} \qquad
A^{op}\to \cC, \quad a\mapsto a\stac B 1_T,
$$
are algebra homomorphisms with commuting ranges in $\cC$, hence $\cC$ is an
$A\sstac k A^{op}$-ring. It remains to check the compatibility between its
$A$-coring and $A\sstac k A^{op}$-ring structures. The Takeuchi property of
the coproduct follows by Definition~\ref{def:ABtor}~(a), its multiplicativity
follows by  Definition~\ref{def:ABtor}~(c) and unitality follows by 
Definition~\ref{def:ABtor}~(d). 
The compatibility of the counit with the multiplication
and unit is obvious. In this way $\cC$ is a right bialgebroid over $A$, and $T$
is its right comodule algebra. By Theorem~\ref{thm.main.1}, $B\subseteq T$ is a
right $\cC$-Galois extension. Since $T$ is a faithfully flat right $A$-module,
this implies that $\cC$ is a right $\times_A$-Hopf algebra by \cite[Lemma
  4.1.21]{Hobst:PhD} (cf. Section \ref{sec.prelims}).

It is proven in a symmetric way that the $B$-coring $\cD$, associated to a
faithfully flat $A$-$B$ torsor $T$, is a left $\times_B$-Hopf algebra and $T$
is its left comodule algebra.

Conversely, let $T$ be a right faithfully flat $A$-ring and a left faithfully
flat right $\cC$-Galois extension of $B$, for a right $\times_A$-Hopf algebra
$\cC$. 
By $B$-$A$ bilinearity, unitality and multiplicativity of the $\cC$-coaction
$\varrho^T$ in $T$, $\varrho^T(ba)=b\sstac A s(a)=\varrho^T(ab)$, where $s$
denotes the source map in $\cC$. Applying the counit we conclude that $B$ and
$A$ are commuting subalgebras of $T$.
Consider the pre-torsor map in Lemma \ref{lemma.tor1}. It satisfies
 Definition~\ref{def:ABtor}~(a) by the Takeuchi property of  $\varrho^T$, and 
 Definition~\ref{def:ABtor}~(b) by its right $B$-linearity.
Definition~\ref{def:ABtor}~(d) follows by the unitality of the right
$\cC$-coaction $\varrho^T$ in $T$, and  Definition~\ref{def:ABtor}~(c) 
follows  by its multiplicativity. 

It is proven in a symmetric way that a left faithfully flat $B$-ring $T$,
which is a right faithfully flat left $\cD$-Galois extension of $A$, for a
left $\times_B$-Hopf algebra $\cD$, is a faithfully flat $A$-$B$ torsor.
This finishes the proof.
\end{proof}

Using the monoidality of the category of comodules of a bialgebroid, one can
find a simpler description of the $\times_A$- and $\times_B$-Hopf algebras,
associated to a faithfully flat $A$-$B$ torsor, than the one in Lemma
\ref{lemma.pure.eq}. Similarly to \cite[Theorem 3.5]{Sch:big}, they turn out to
be coinvariants of diagonal comodules. 
\begin{lemma}\label{lem:TTcoC}
Let $T$ be a faithfully flat $A$-$B$ torsor with associated right
$\times_A$-Hopf algebra $\cC$ and left $\times_B$-Hopf algebra $\cD$. 
\begin{zlist}
\item View $T\sstac A T$ as a right $\cC$-comodule with the diagonal coaction 
\begin{equation}\label{eq:Cdiag}
T\stac A T\to T\stac A T\stac A \cC,\qquad u\stac A v\mapsto 
u\sut 1 \stac A v\sut 1 \stac A v\sut 2 u\sut 2 \stac B u\sut 3 v\sut 3.
\end{equation}
Then $\cD = (T\sstac AT)^{co\cC}$.
\item View $T\sstac B T$ as a left $\cD$-comodule with the diagonal coaction 
\begin{equation}\label{eq:Ddiag}
T\stac B T\to \cD \stac B T\stac B T,\qquad u\stac B v\mapsto 
u\sut 1 v\sut 1 \stac A v\sut 2 u\sut 2 \stac B u\sut 3 \stac B v\sut 3.
\end{equation}
Then $\cC ={}^{co\cD}\!(T\sstac B T)$.
\end{zlist}
\end{lemma}
\begin{proof}
(1) Recall from the proof of Theorem \ref{thm.main.1} that an element $\sum_i
    u_i\sstac A v_i\in T\sstac A T$ belongs to $\cD$ if and only if 
\begin{equation}\label{eq:TTcoC1}
 \sum_i u_i\sut 1\stac A u_i\sut 2\stac B u_i\sut 3 v_i
=\sum_i u_i\stac A v_i\stac B 1_T.
\end{equation}
Application of the bijective canonical map $T\sstac B\, T\to T\sstac A\, \cC$,
$t\sstac B\, t'\mapsto t t'\sut 1 \sstac A t'\sut 2\sstac B t'\sut 3$ to
the last two factors of \eqref{eq:TTcoC1} yields the equivalent condition
\begin{equation}\label{eq:TTcoC2}
\sum_i u_i\sut 1 \stac A u_i\sut 2 (u_i\sut 3 v_i)\sut 1\stac A 
(u_i\sut 3 v_i)\sut 2\stac B (u_i\sut 3 v_i)\sut 3 =
\sum_i u_i\stac A v_i 1_T\sut 1  \stac A 1_T\sut 2 \stac B 1_T \sut 3.
\end{equation}
By  Definition~\ref{def:ABtor}~(c) and (d),  and Definition~\ref{def.pre.tor}
(c) and (b), \eqref{eq:TTcoC2} is equivalent to 
\begin{equation}\label{eq:TTcoC3}
\sum_i u_i \sut 1 \stac A v_i \sut 1 \stac A  v_i \sut 2 u_i \sut 2 \stac B
u_i\sut 3 v_i\sut 3= \sum_i u_i\stac A v_i \stac A 1_T \stac B 1_T.
\end{equation}
Equation \eqref{eq:TTcoC3} expresses the property that $\sum_i u_i\sstac A
v_i\in T\sstac A T$ is coinvariant with respect to the coaction
\eqref{eq:Cdiag}. This proves claim (1). Part (2) is proven by a symmetrical
reasoning.
\end{proof}

Consider a right bialgebroid $\cC=(\cC,s,t,\Delta_\cC,\varepsilon_\cC)$ over a
$k$-algebra $A$ and
a $k$-algebra $B$. Write $B^e = B\sstac k B^{op}$ for the enveloping algebra
of $B$ and view any $B$-$B$ bimodule as a left $B^e$-module. 
Let $T$ be a $B^e$-$\cC$ bicomodule.
Then for every left $\cC$-comodule $M$, the cotensor product $T\Box_\cC
M$ inherits a left $B^e$-module structure of $T$, i.e.\ $T$ {\em induces a
  functor} 
$U_T :=T\Box_\cC\, \bullet:{}^\cC \cM\to {}_B \cM_B$. Our next task is a
study of this functor. Our line of reasoning follows ideas in
\cite{Ulb:Galfun}, although we have to face complications caused by the
fact that the algebras $A$ and $B$ might be different from $k$.
Recall from \cite[39.3]{BrzWis:cor} that if $U_T$
preserves colimits, then $(T\Box_\cC M)\sstac R N \cong T\Box_\cC (M\sstac R
N)$ canonically, for any algebra $R$, $\cC$-$R$ bicomodule $M$ and left
$R$-module $N$. 

\begin{theorem}\label{thm:lax}
Let $B$ be an algebra, $\cC=(\cC,s,t,\Delta_\cC,\varepsilon_\cC)$ a right
$A$-bialgebroid and $T$ a $B^e$-$\cC$ bicomodule. Let
$U_T :=T\Box_\cC\, \bullet:{}^\cC \cM\to {}_B \cM_B$.
\begin{zlist}
\item Every $\cC$-comodule algebra structure in $T$ which makes $T$ a $B$-ring
  determines a lax monoidal structure of $U_T$.
\item If $U_T$ preserves colimits, then every lax monoidal structure of $U_T$
  determines a $\cC$-comodule algebra structure in $T$, such that $T$ is a
  $B$-ring.  
\item If $U_T$ preserves colimits, then the constructions in parts (1) and (2)
  are mutual inverses.  
\end{zlist}
\end{theorem}
\begin{proof}
(1) The proof consists of a construction of natural homomorphisms $\xi_0:B\to
    T\Box_\cC A$ and $\xi_{\bullet ,\bullet}:(T\Box_\cC\, \bullet)\sstac B 
    (T\Box_\cC\, \bullet)\to T\Box_\cC (\bullet \sstac {A^{op}}\, \bullet)$,
    making $U_T$ a monoidal functor. Consider the maps 
\begin{eqnarray}
&&\xi_0(b) := 1_Tb\stac A 1_A,\label{eq:xi0}\\
&&\xi_{M,M'}(\sum_i (u_i\stac A m_i)\stac B (u'_i\stac A m'_i)) := \sum_i
u_i u'_i \stac A (m_i\stac {A^{op}} m'_i),\label{eq:xi}
\end{eqnarray}
for any left $\cC$-comodules $M$ and $M'$. Since the cotensor product
$T\Box_\cC \ M'$ is contained in the centraliser of $A$ in the obvious $A$-$A$
bimodule $T\sstac A M'$, it is 
easy to see that $\xi_{M,M'}$ is well defined. By the left $B$-linearity and
unitality of the coaction $\varrho^T$, and unitality of the target map $t$,
the range of $\xi_0$ is in the required cotensor product. By the
multiplicativity of $\varrho^T$, also the range of $\xi_{M,M'}$ is in the
appropriate cotensor product. Obviously, both $\xi_0$ and $\xi_{M,M'}$ are
$B$-$B$ bilinear. Naturality of $\xi_{\bullet,\bullet}$ (in both arguments)
follows easily by its explicit form. The hexagon identity follows by the
associativity of the multiplication in $T$ and the triangle identities follow
by its unitality.

(2) The proof consists of a construction of right $\cC$-colinear
    multiplication and unit maps, equipping $T$ with $A$- and $B$-ring
    structures.  
Note first that the lax monoidal functor $U_T$ maps the algebra $\cC$ in
${}^\cC \cM$ (with unit map $t$ and coaction $\Delta_\cC$) to an algebra
$T\Box_\cC\, \cC \cong T$ in ${}_B \cM_B$, with 
structure maps 
\begin{eqnarray*}
&&\hspace{-.8cm}\xymatrix{
\mu_T:\ \big(T\stac B T \ar[rr]^{\varrho^T\sstac B \varrho^T}&&
( T\Box_\cC \, \cC)\sstac B ( T\Box_\cC \, \cC) \ar[r]^{\ \
  \quad\xi_{\cC,\cC}}& 
T\cstac \cC \, (\cC\stac{A^{op}} \, \cC)\ar[r]^{\ \quad T\Box_\cC \mu_\cC}&
T\cstac \cC \cC \ar[r]^{\quad T\sstac A \varepsilon_\cC}&T\big),
}
\\
&&\hspace{-.8cm}\xymatrix{
\eta_T:\ \big( B\ar[r]^{\xi_0}&
T\cstac \cC A \ar[r]^{T\Box_\cC t}&
T\cstac \cC  \cC \ar[r]^{\quad T\sstac A \varepsilon_\cC}& 
T\big).
}
\end{eqnarray*}
This proves that $T$ is a $B$-ring or, equivalently,
that $T$ is a $k$-algebra with unit $1_T:= \eta_T(1_B)$ and $\eta_T:B\to
T$ is a homomorphism of  $k$-algebras. 
We make no notational difference between the multiplication maps in $T$ as a 
$k$-algebra or $B$-ring. 
In all cases multiplication will be denoted by juxtaposition. It should be
clear from the context, which structure is meant.

In order to show that $T$ is also a $\cC$-comodule algebra, we investigate 
properties of  the map $\xi_{\cC,\cC}$.
Considering $\cC\sstac A\cC$ as a left $\cC$-comodule via the regular comodule
structure of the first factor, the coproduct $\Delta_\cC$ is left
$\cC$-colinear. Hence the naturality of $\xi_{\bullet,\bullet}$ implies the
  identities
\begin{eqnarray}\label{eq:col1}
\big( T\cstac \cC (\cC \stac {A^{op}} \Delta_\cC)\big)\circ \xi_{\cC,\cC}&=&
\xi_{\cC,\cC\sstac A \cC}\circ \big( (T\cstac \cC \cC )\stac B (T\cstac \cC
\Delta_\cC ) \big),\\\label{eq:col2}
\big( T\cstac \cC (\Delta_\cC \stac {A^{op}} \cC )\big)\circ \xi_{\cC,\cC}&=&
\xi_{\cC\sstac A \cC, \cC}\circ \big( (T\cstac \cC \Delta_\cC)\stac B (T\cstac
\cC \cC ) \big).
\end{eqnarray}
For any left $A$-module $N$, $\cC \sstac A N$ is a left
$\cC$-comodule via the first factor and the map $\cC\to \cC \sstac A N$,
$c\mapsto c\sstac A n$ is left $\cC$-colinear, for any $n\in N$. Hence the
naturality of $\xi_{\bullet,\bullet}$ implies that, for any left
$\cC$-comodule $M$, the following diagrams are commutative.
\begin{equation}\label{eq:col3}
\xymatrix{
(T\cstac \cC M)\stac B (T\cstac \cC \cC)\stac A N
\ar[r]^{\cong}
\ar[d]_{\xi_{M,\cC}\stac A N}&
(T\cstac \cC M)\stac B \big(T\cstac \cC (\cC\stac A N)\big) 
\ar[d]^{\xi_{M,\cC\sstac A N}}\\
\big(T \cstac \cC(M\stac {A^{op}} \cC)\big)\stac A N \ar[r]^{\cong}&
T\cstac \cC \big(M\stac {A^{op}}\cC \stac A N\big)
}
\end{equation}
and
\begin{equation}\label{eq:col4}
\xymatrix{
\big((T\cstac \cC \cC_\bullet)\stac B (T\cstac \cC M)\big)\stac A N
\ar[r]^{\cong} 
\ar[d]^{\xi_{\cC_\bullet,M}\stac A N}&
\big(T\cstac \cC (\cC \stac A N)\big)\stac B (T\cstac \cC M) 
\ar[d]^{\xi_{\cC\sstac A N,M}}\\
\big(T \cstac \cC(\cC_\bullet  \stac {A^{op}} M)\big)\stac A N
\ar[r]^{\cong} &
T\cstac \cC \big((\cC_\bullet\stac A N) \stac {A^{op}} M\big )
}
\end{equation}
The horizontal arrows are isomorphisms by the assumption that the functor
$T\Box_\cC\, \bullet $ preserves colimits.
Put $M=N=\cC$ in \eqref{eq:col3} and \eqref{eq:col4}. Combining the resulting
equalities with \eqref{eq:col1} and \eqref{eq:col2}, we conclude that, for
$\sum_i (u_i\sstac A c_i)\sstac B (u'_i\sstac A c'_i)\in (T\Box_\cC\,
\cC)\sstac B (T\Box_\cC\, \cC)$, 
\begin{equation}\label{eq:xi_CC}
\sum_i u_iu'_i\stac A c_i \stac {A^{op}} c'_i = \xi_{\cC,\cC}\big(\sum_i
(u_i\stac A c_i)\stac B (u'_i\stac A c'_i) \big),
\end{equation}
proving that the left hand side is well defined.
(A more detailed version of a more general computation will be presented in
the proof of part (3).)

In order to check that the coaction $\varrho^T$ is unital, introduce
the notation $\sum_k u_k\sstac A a_k  := \xi_0(1_B)\in T\Box_\cC A$. By
the definition of the unit in the $B$-ring $T$, $\sum_k u_k a_k=1_T$.  With
this identity at hand,  
\begin{equation}\label{eq:rho(1)}
\varrho^T(1_T)= \sum_k u_k\su 0 \stac A u_k\su 1 s(a_k)=
\sum_k u_k \stac A t(a_k)=\sum_k u_k a_k \stac A 1_\cC =1_T\stac A 1_\cC,
\end{equation}
where $\roT(u) = u\su 0\sstac{A}u\su 1$ is the Sweedler notation for the right 
$\cC$-coaction.
The first equality follows by the right $A$-linearity of $\varrho^T$.
In the second equality we used the fact that $\sum_k u_k\sstac A a_k=
\xi_0(1_B)$ belongs to the cotensor product. 

Note that by the right $A$-linearity of the coaction $\varrho^T$, naturality
of $\xi_{\bullet,\bullet}$ (cf. \eqref{eq:col3} for $N=A$), and right
$A$-linearity of the multiplication 
$\mu_\cC$ and counit $\varepsilon_\cC$ in $\cC$, the multiplication map
$\mu_T$ is right $A$-linear. Hence, for all $a,a'\in A$,
$$
\mu_T(1_T a\stac B 1_T a')=\mu_T(1_T a\stac B 1_T) a'=(1_T a)a'=1_T(aa'),
$$
where the second equality follows by the unitality of $\mu_T$. Therefore the
map
$$
A\to T,\qquad a\mapsto 1_T a 
$$
is an algebra homomorphism and $ta=t(1_T a)$. Furthermore, by right
$A$-linearity and unitality of $\varrho^T$ and \eqref{eq:xi_CC}, 
$$
(1_Ta)t=(T\stac A \eC\circ \mu_\cC)\circ \xi_{\cC,\cC}((1_T\stac A s(a))\stac
B  (t\su 0\stac A t\su 1))=t\su 0\eC(s(a)t\su 1)=at.
$$
This proves that $T$ (with the $A$-actions determined by its  $\cC$-comodule
structure) is an $A$-ring.  
Furthermore, \eqref{eq:xi_CC} implies, for $u,v\in T$,
\begin{eqnarray*}
u\su 0 v\su 0\stac A u\su 1 v\su 1
&=& \big((T\cstac\cC\mu_\cC)\circ \xi_{\cC,\cC}\big)\big((u\su 0\stac A u\su
1)\stac B (v\su 0\stac A v\su 1)\big)\\
&=& \big(\varrho^T\circ (T\stac A \varepsilon_\cC\circ \mu_\cC)\circ
\xi_{\cC,\cC}\circ (\varrho^T\stac B \varrho^T)\big)(u\stac B v)
=(uv)\su 0\stac A (uv)\su 1.
\end{eqnarray*}
This proves the right $\cC$-colinearity of the multiplication map in the
$A$-ring $T$.
The right $\cC$-colinearity of the unit map is equivalent to the unitality of
the coaction $\varrho^T:T\to T\times_A \cC$, proven in \eqref{eq:rho(1)}. 
Hence $T$ is a right $\cC$-comodule algebra, as stated. This completes the
proof of part (2). 

(3) It is straightforward to see that starting with a right $\cC$-comodule
    algebra $T$, which is also a $B$-ring, and applying first the
    construction in part (1) to it and then the one in part (2) to the result,
    one recovers the original algebra structure in $T$.

Conversely, starting with a lax monoidal structure on $U_T$, with coherence
    natural homomorphisms $\xi_0$ and $\xi_{\bullet,\bullet}$, applying first
    the construction in part (2) to it and then the one in part (1) to the
    result, one recovers the original natural homomorphism $\xi_0$ and a 
    natural homomorphism 
\begin{eqnarray}\label{eq:newxi}
(T\cstac \cC M)\stac B (T\cstac \cC M') &\to& T\cstac \cC (M\stac {A^{op}} M'),
\\
\sum_i (u_i \stac A m_i)\stac B (u'_i\stac A m'_i) &\mapsto&\nonumber\\
\sum_i \big( (T\stac A \varepsilon_\cC \circ \mu_\cC)&\hspace{-.8cm}\circ&
\hspace{-.8cm}\xi_{\cC,\cC}\big) 
\big((u_i\su 0\stac A u_i \su 1)\stac B (u'_i \su 0 \stac A u'_i \su 1)\big)
\stac A (m_i \stac {A^{op}} m'_i)\equiv \nonumber\\
\sum_i \big( (T\stac A \varepsilon_\cC \circ
\mu_\cC)&\hspace{-.8cm}\circ&\hspace{-.8cm} \xi_{\cC,\cC}\big) 
\big((u_i\stac A m_i\su {-1})\stac B (u'_i\stac A m'_i \su{-1})\big)\stac A
(m_i\su 0\stac {A^{op}} m'_i \su 0), \nonumber
\end{eqnarray}
where, for all $m\in M$, ${}^M\!\varrho(m) = m \su{-1}\sstac A m\su 0$ is the
notation 
for the left coaction.
The claim is proven by showing that the map \eqref{eq:newxi} is
equal to $\xi_{M,M'}$. This can be seen by similar arguments to those used in
proving part (2). By the left $\cC$-colinearity of the left $\cC$-coaction
${}^M \varrho:M\to \cC\sstac A M$ (where the codomain is a left $\cC$-comodule
via the first factor) and naturality of $\xi_{\bullet,\bullet}$, 
\begin{eqnarray}
\big(T\cstac \cC(M\stac {A^{op}} {}^{M'}\varrho)\big)\circ \xi_{M,M'}&=&
\xi_{M, \cC\sstac A M'}\circ \big( (T\cstac \cC M)\stac B (T\cstac \cC
{}^{M'}\varrho)\big) \label{eq:xinat2}\\
\big(T\cstac \cC({}^M\varrho\stac {A^{op}} {M'})\big)\circ \xi_{M,M'}&=&
\xi_{\cC\sstac A M, M'}\circ \big( (T\cstac \cC {}^M \varrho)\stac B 
(T\cstac \cC {M'})\big)\label{eq:xinat1}.
\end{eqnarray}
Introduce the notation $\xi_{M,M'}\big(\sum_i (u_i\sstac A m_i)\sstac B
(u'_i\sstac A m'_i)\big)=\colon \sum_j t_j \sstac A (n_j\sstac {A^{op}} n'_j)$,
for any fixed element $\sum_i (u_i\sstac A m_i)\sstac B (u'_i\sstac A m'_i)\in
(T\Box_\cC M)\sstac B (T\Box_\cC M')$. The identities \eqref{eq:col4},
\eqref{eq:xinat1}, \eqref{eq:col3} and \eqref{eq:xinat2} together imply the 
identity
\begin{eqnarray}\label{eq:bigdiagr}
&&\sum_i \big(\xi_{\cC,\cC}\big((u_i\stac A m_i\su{-1})\stac B (u'_i\stac A
m'_i\su{-1})\big)\stac A m_i\su 0\big)\stac A m'\su 0=\nonumber\\
&& \sum_j\big(\big(t_j \stac A (n_j\su{-1}\stac {A^{op}}
n'_j\su{-1})\big)\stac A n_j\su 0\big)\stac A n'_j\su 0
\end{eqnarray}
in $\big(\big(T\sstac A (\cC\sstac {A^{op}}\cC)\big)\sstac A M\big)\times_A
  M'$. Recall that in the tensor product $\cC\sstac {A^{op}}\cC$ both module
  structures are given by the target map. In the tensor product 
$T\sstac A (\cC\sstac {A^{op}}\cC)$ the left $A$-module structure of
  $\cC\sstac {A^{op}}\cC$, given by right multiplication by the target map in
  the second factor is used. In $\big(T\sstac A (\cC\sstac
  {A^{op}}\cC)\big)\sstac A M$ the first tensorand $T\sstac A (\cC\sstac
  {A^{op}}\cC)$ is understood to be a right $A$-module via right
  multiplication by the source map in the middle factor $\cC$. Finally, the
  Takeuchi product $\big(\big(T\sstac A (\cC\sstac {A^{op}}\cC)\big)\sstac A
  M\big)\times_A M'$ is the center of the $A$-$A$ bimodule $\big(\big(T\sstac
  A (\cC\sstac {A^{op}}\cC)\big)\sstac A M\big)\sstac A M'$, where
  $(\big(T\sstac A (\cC\sstac {A^{op}}\cC)\big)\sstac A M$ is an $A$-$A$
  bimodule via left and right multiplications by the source map in
  the third factor $\cC$. Consider the map
\begin{eqnarray}\label{eq:rebr}
&\big(\big(T\stac A (\cC\stac {A^{op}}\cC)\big)\stac A M\big)\times_A M'
\to &T\stac A (\cC\stac {A\otimes A^{op}} \cC)\stac A (M\stac {A^{op}} M'),
\nonumber\\
& \big(\big(t\stac A (c\stac {A^{op}} c')\big)\stac A m\big) \stac A m'\mapsto
& t\stac A (c\stac {A\otimes A^{op}} c')\stac A (m \stac {A^{op}} m').
\end{eqnarray} 
In its codomain the tensor product $\cC\sstac {A\otimes A^{op}} \cC$
corresponds to the $A\otimes A^{op}$-ring structure of $\cC$, via the source
and target maps. In the tensor product $T\sstac A (\cC\sstac {A^{op}}\cC)$ the
left $A$-module structure of $\cC\sstac {A^{op}}\cC$ is given by right 
multiplication by the target map in the second factor, as in the domain. In
the tensor product $(\cC\sstac {A\otimes A^{op}} \cC)\sstac A (M\sstac {A^{op}}
M')$ the first factor $\cC\sstac {A\otimes A^{op}} \cC$ is understood to be a
right $A$-module via right multiplication by the source map in the second
factor, and $M\sstac {A^{op}} M'$ is a left $A$-module via the left $A$-module
structure of $M'$. Apply the map \eqref{eq:rebr} to both sides of
\eqref{eq:bigdiagr} to conclude that the map \eqref{eq:newxi} takes $\sum_i
(u_i\sstac A m_i)\sstac B (u'_i\sstac A m'_i)\in (T\Box_\cC M)\sstac B
(T\Box_\cC M')$ to
\begin{eqnarray*}
&&\sum_j (T\stac A \varepsilon_\cC\circ \mu_\cC)(t_j \stac A (n_j
  \su{-1}\stac{A\otimes A^{op}} n'_j\su{-1}))\stac A (n_j\su 0\stac {A^{op}}
  n'_j\su 0)=\\
&&\sum_j t_j \stac A \varepsilon_\cC(n_j \su {-1} n'_j\su{-1}))(n_j\su 0\stac
  {A^{op}}   n'_j\su 0)=
\sum_j t_j \stac A \varepsilon_\cC\big((n_j\stac{A^{op}}   n'_j)\su
  {-1}\big)(n_j\stac{A^{op}}   n'_j)\su 0=\\
&&\sum_j t_j \stac A (n_j\stac{A^{op}}   n'_j)
=\xi_{M,M'}\big(u_i\sstac A m_i)\sstac B (u'_i\sstac A m'_i)\big).
\end{eqnarray*}
This completes the proof.
\end{proof}

\begin{lemma}\label{lem:iso}
Let $\cC=(\cC,s,t,\Delta_\cC,\varepsilon_\cC)$ be a right $\times_A$-Hopf
algebra. For any left $A$-modules $N$ and $M$, and right $\cC$-comodule $T$,
the map  
\begin{equation}\label{eq:iso}
T\cstac \cC\big( (\varepsilon_\cC\stac A N)\stac {A^{op}} (\cC\stac A M)\big):
T\cstac \cC\big( (\cC\stac A N)\stac {A^{op}} (\cC\stac A M)\big)\to
(T\stac A \cC_{\bullet} \stac A N)\stac A M,
\end{equation}
is a bijection. In the domain of \eqref{eq:iso}, $\cC$ is a right $A$-module
through  the right multiplication by the source map,
$\cC\sstac A M$ and $\cC\sstac A N$ are left $\cC$-comodules via the regular
comodule structure of their first factors, 
and $(\cC\sstac A N)\sstac {A^{op}} (\cC\sstac A M)$ is a left
$\cC$-comodule via the diagonal coaction. In the codomain, $\cC$ in
$\cC\sstac A N$ is  the right $A$-module 
through the left multiplication by the target map, while in  $\cC\sstac A M$
it is a right $A$-module through the right multiplication by the source map. In
$T\sstac A \cC$ the left $A$-module structure of $\cC$ is understood via right
multiplication by the target map.
\end{lemma}
\begin{proof}
The map \eqref{eq:iso} is well defined by the $A$-$A$ bilinearity of
$\varepsilon_\cC$. Its bijectivity is proven by constructing the inverse in
terms 
of the inverse of the Galois map $\theta: \cC\sstac {A^{op}} \cC\to \cC\sstac
A \cC$, $c\sstac A c'\mapsto c\Delta_\cC(c')$. 
 Consider the map
\begin{eqnarray}\label{eq:inv}
\\
&&\big((T\stac A \mu_\cC \stac{A^{op}} \cC)\circ (\varrho^T\stac A
\theta^{-1}(1_\cC\stac A \, \bullet))\stac A N\big)\stac A M:\nonumber\\
&&(T\stac A \cC_\bullet\stac A N)\stac A M\to 
\big((T\stac A \cC_\bullet\stac {A^{op}} \cC_\circ )\stac A {}_\bullet
N\big)\stac A {}_\circ M
\cong T\stac A \big( (\cC\stac A N)\stac {A^{op}} (\cC \stac A M)\big).
\nonumber
\end{eqnarray}
It is well defined by the multiplicativity of the translation map
$\theta^{-1}(1_\cC\sstac A \, \bullet)$ and \eqref{eq:thetaA}.
Application of the pentagon identity \eqref{eq:pent} in the form 
$$
\theta_{13}^{-1} \circ (\cC \stac A \theta^{-1})\circ (\theta\stac A \cC)=
(\theta\stac {A^{op}} \cC)\circ (\cC \stac {A^{op}} \theta^{-1})
$$
yields that
the range of \eqref{eq:inv} is 
in the cotensor product $T\Box_ \cC  \big( (\cC\sstac A N)\sstac {A^{op}} (\cC
\sstac A M)\big)$. We leave it to the reader to check that the maps
\eqref{eq:iso} and \eqref{eq:inv} are mutual inverses.
\end{proof}

\begin{theorem}\label{thm:strmon}
Let $\cC=(\cC,s,t,\Delta_\cC,\varepsilon_\cC)$ be a right $\times_A$-Hopf
algebra. Let $T$ be a right  
$\cC$-comodule algebra and $B$ a subalgebra of $T^{co\cC}$. Denote by   $U_T
:= T\Box_\cC\, \bullet$ the induced lax monoidal  functor ${}^\cC\!\cM\to
{}_B\cM_B$. 
\begin{zlist}
\item If $U_T$ is a monoidal functor, then $B\subseteq T$ is a right
  $\cC$-Galois extension.
\item If $B\subseteq T$ is a right faithfully flat right $\cC$-Galois
  extension, then $U_T$ is a monoidal functor.
\end{zlist}
\end{theorem}
\begin{proof}
(1) Via the embedding $T\Box_\cC A \hookrightarrow T\sstac A A \cong T$, the
coinvariants of the right $\cC$-comodule $T$ are identified with the elements
of $T\Box_\cC A$. By the monoidality of $U_T$, $\xi_0:B\to T\Box_\cC A$
is an isomorphism. Hence $B$ is equal to $T^{co\cC}$. 

In light of the form \eqref{eq:xi} of the natural homomorphism
$\xi_{\bullet,\bullet}:(T\Box_\cC \, \bullet)\sstac B (T\Box_\cC \,
\bullet)\to T\Box_\cC \ (\bullet\sstac {A^{op}}\, \bullet)$, the canonical 
map $T\sstac B T\to T\sstac A \cC$ can be written as a composite map,
\begin{equation}\label{eq:can}
\xymatrix{
T\stac B T\ar[rr]^{\varrho^T\sstac B \varrho^T\quad}&&
(T\cstac \cC \cC)\stac B (T\cstac \cC \cC)\ar[r]^{\quad\xi_{\cC,\cC}}&
T\cstac \cC (\cC \stac {A^{op}} \cC)\ar[rr]^{\quad T\Box_\cC
  (\varepsilon_\cC\sstac   {A^{op}} \cC)}&&T\stac A \cC.
}
\end{equation}
Since $\varrho^T:T\to T\Box_\cC\, \cC$ is an isomorphism, so is the first
arrow in \eqref{eq:can}. The middle arrow is an isomorphism by the monoidality
of $U_T$. The last arrow is an isomorphism by Lemma~\ref{lem:iso}. This
proves that the canonical map \eqref{eq:can} is bijective, that is,
$B\subseteq T$ is a right $\cC$-Galois extension.

(2) We need to show that the natural homomorphisms \eqref{eq:xi0} and
    \eqref{eq:xi} are isomorphisms. Since the canonical map $\can:T\sstac B
    T\to T\sstac A \cC$ is an isomorphism by assumption, so is the map
$$
\xymatrix{
\gamma_M:\ T\stac B(T\cstac \cC M)\cong (T\stac B T)\cstac \cC M
\ar[rr]^{\quad \can \Box_\cC M}&&(T\stac A \cC)\cstac \cC M\cong T\stac A
M\ ,} 
$$
mapping $\sum_i u_i \sstac B v_i \sstac A m_i$ to $\sum_i u_iv_i\sstac A m_i$,
for any left $\cC$-comodule $M$. Since $T\Box_\cC\, M$ is contained in the
centraliser of $A$ in the obvious $A$-$A$ bimodule $T\sstac A M$, the map
$\gamma_M$ is right $A$-linear with respect to the
right $A$-module structure of $T\sstac B (T\Box_\cC \, M)$, via its first
factor. A straightforward computation shows the
commutativity of the following diagram, for any left $\cC$-comodules $M$ and
$M'$, 
$$
\xymatrix{
T\stac B (T\cstac \cC M)\stac B(T\cstac\cC M')\ar[rrr]^{\gamma_M \sstac
  B(T\Box_\cC M')}\ar[d]^{T\stac B \xi_{M,M'}}&&&
\stackrel{\displaystyle (T\stac A M)\stac B (T\cstac\cC M')}{\cong (T_\bullet
  \stac B(T\cstac   \cC M'))\stac A M}\ar[d]^{\gamma_{M'}\stac A M}\\
T\stac B \big(T\cstac \cC (M\stac {A^{op}} M')\big)\ar[rrr]^{\gamma_{M\sstac
    {A^{op}} M'}}&&&
\stackrel{\displaystyle T\stac A M'\stac A M}{\cong T\stac A (M\stac {A^{op}}
  {}_\bullet   M')}
}
$$
Therefore $T\sstac B \xi_{M,M'}$ is an isomorphism. Since $T$ is a faithfully
flat right $B$-module by assumption, this proves that $\xi_{M,M'}$ is an
isomorphism, for any left $\cC$-comodules $M$ and $M'$. By one of the triangle
identities, $\xi_{\cC,A}\circ (\varrho^T\sstac B (T\Box_\cC A))\circ (T\sstac
B \xi_0)$ is an isomorphism. Hence so is $T\sstac B \xi_0$, and, by the
faithful flatness of the right $B$-module $T$, also $\xi_0$. This completes
the proof.
\end{proof}

In particular Theorem~\ref{thm:strmon} implies that a faithfully flat 
$A$-$B$ torsor which is also faithfully flat as a right $B$-module   induces
a $k$-linear monoidal functor ${}^\cC\cM\to {}_B \cM_B$. 
In the case when $A$ and $B$ are both equal to the ground ring $k$ (and hence
the faithful flatness assumptions made on $T$ imply properties (iv) in Remark
\ref{rem:ff}), the induced functor is known to be a `fibre functor'. That is,
it is faithful and preserves colimits and kernels. Actually it is not hard to
see that also for arbitrary $k$-algebras $A$ and $B$, properties (iv) in Remark
\ref{rem:ff} imply these properties of the induced functor. Furthermore, in
the case when $A$ and $B$ are trivial, every fibre functor is known to be
induced by a faithfully flat torsor. In our general setting it follows by
Theorems \ref{thm:lax} and \ref{thm:strmon} that, for a right $\times_A$-Hopf
algebra $\cC$, an algebra $B$, and a $k$-linear faithful monoidal functor
$U:{}^\cC\cM\to{}_B \cM_B$, preserving colimits and kernels, $U(\cC)$ is a
right $\cC$-Galois extension of $B$. However, we were not able
to derive properties (iv) in Remark \ref{rem:ff} of $U(\cC)$.  

Finally, we show that the functor induced by 
a torsor satisfying property (iii) in  Remark \ref{rem:ff} is a monoidal
equivalence. 

\begin{theorem}\label{thm:tor-moneq}
Let  $T$ be
an $A$-$B$ torsor which obeys properties (iii) in Remark~\ref{rem:ff}. 
  Then the functor 
  $T\Box_\cC\,\bullet:{}^\cC\cM\to {}^\cD \cM$ is a 
  $k$-linear monoidal equivalence between the categories of comodules for the
  associated right $\times_A$-Hopf algebra $\cC$ and left 
  $\times_B$-Hopf algebra $\cD$.
\end{theorem}
\begin{proof}
The functor $T\Box_\cC\, \bullet:{}^\cC\cM\to {}^\cD\cM$ is a $k$-linear
    equivalence by 
    Corollary~\ref{cor:prtoreq}. The functor $T\Box_\cC\,\bullet:{}^\cC\cM\to
    {}_B \cM_B$ is monoidal by Theorem \ref{thm:strmon} (2). Recall the
    existence of a strict monoidal forgetful functor ${}^\cD\cM \to {}_B
    \cM_B$. It is
    straightforward to see that the left $\cD$-comodule algebra structure of
    $T$ implies the left $\cD$-colinearity of the coherence natural
    isomorphisms \eqref{eq:xi0} and \eqref{eq:xi}. This completes the
    proof.
 \end{proof}

Note that in the case when the algebras $A$ and $B$ are equal to the ground
ring $k$, also a converse of Theorem \ref{thm:tor-moneq} holds: every
$k$-linear monoidal
equivalence between comodule categories of flat Hopf algebras is induced by a
faithfully flat torsor \cite[Corollary 5.7]{Sch:big}. In our general setting,
however, by 
the failure of properties (iii) in Remark \ref{rem:ff} of $F(\cC)$, for a
$k$-linear monoidal equivalence functor $F:{}^\cC \cM \to {}^\cD \cM$, we were
not able to prove an analogous result.
\appendix
\section{Cleft extensions by Hopf algebroids}
A Hopf algebroid consists of a compatible pair of a right and a left
bialgebroid -- on the same total algebra over a commutative ring $k$ but over
anti-isomorphic base $k$-algebras $A$ and $L$. In addition there is an
antipode map. References about Hopf algebroids are \cite{BohmSzl:hgdax} (in the
case when the antipode is bijective) and \cite{Bohm:hgdint} (in general). As a
consequence of the compatibility of the two involved bialgebroid structures,
the categories of their (right, say) comodules are isomorphic monoidal
categories. The isomorphism is compatible with the forgetful functors to the
category of $k$-modules. Its explicit form was given in \cite[Theorem
  2.2]{BohBrz:cleft}.  

Since in a Hopf algebroid there are two bialgebroid (and hence coring)
structures present, it is helpful to use two versions of Sweedler's index
notation.  
Upper indices are used to denote the components of the coproduct 
in the right bialgebroid, i.e.\ we write $h\mapsto h^{(1)}\otimes_A h^{(2)}$. 
Lower indices,   $h\mapsto h_{(1)}\otimes_L h_{(2)}$, denote the coproduct of
left bialgebroid. Similar upper/lower index conventions are used for right
coactions of the right/left bialgebroid (related by the isomorphism in
\cite[Theorem 2.2]{BohBrz:cleft}), i.e.\ notations $t\mapsto t\su 0\otimes_A
t\su 1$ / $t\mapsto t_{[0]}\otimes_L t_{[1]}$ are used.

Cleft extensions by Hopf algebroids were introduced in \cite{BohBrz:cleft}. 
In \cite[Theorem 3.12]{BohBrz:cleft} these extensions were characterised as
Galois extensions by the constituent right bialgebroid that satisfy normal
basis property with respect to the constituent left bialgebroid.
\smallskip

Consider a Hopf algebroid $\hH$, with constituent right bialgebroid over an
algebra $A$ and left bialgebroid over $L$.
Let $T$ be a right faithfully flat $A$-ring and a left faithfully flat right
{\em cleft} extension of $B$ by $\hH$.
In this situation, by Theorem~\ref{thm.main.Hobst} and 
Lemma~\ref{lem:TTcoC}, there exists a canonical left $\times_B$-Hopf
algebra structure in $(T\sstac A T)^{co\hH}$. 
The aim of this section is to determine it explicitly.
This results
in new examples of bialgebroids which generalise the extended Hopf algebra 
of Connes and Moscovici introduced in \cite{ConMos:Ran} in the
context of transverse geometry. We start
with the following generalisation of a construction in \cite{Kad:Pseudo} (see
also \cite{PaVOys:K-CMiso}). For definitions of an algebra measured by a
bialgebroid, an (invertible) 2-cocycle on a bialgebroid, a cocycle twisted
module of a bialgebroid, a crossed product by
a bialgebroid and a cleft extension by a Hopf algebroid, we refer to
\cite{BohBrz:cleft}. 
\begin{proposition}\label{prop:twbgd}
Let $\hH$ be a left $\times_L$-Hopf algebra, and $B$ an $L$-ring measured by
$\hH$, with measuring $\cdot$. Let $\sigma$ be a $B$-valued 2-cocycle on
$\hH$, that makes $B$ a 
$\sigma$-twisted $\hH$-module. Consider the $k$-module $\cD := B\sstac L
(B\sstac L {}_\bullet \hH)$, where the $L$-module tensor product in the
parenthesis is understood with respect to the left $L$-module structure of
$\hH$, given through the right multiplication by the target map, and the
resulting 
tensor product is meant to be a left $L$-module via left multiplication by the
source map in the second factor $\hH$. It has the following structures.
\begin{zlist}
\item $\cD$ is a $B\sstac k B^{op}$-ring, with source and target maps
$$
s_\cD(b) := b\stac L 1_B \stac L 1_\hH\qquad \textrm{and}\qquad 
t_\cD(b) := 1_B \stac L b\stac L 1_\hH,
$$
respectively, and multiplication 
\begin{eqnarray*}
(b\stac L b'\stac L h)&&\hspace{-.5cm}(c\stac L c'\stac L k) :=\\
&&b(h\sut 1 \sw 1\cdot c)\sigma(h\sut 1 \sw 2,k\sut 1 \sw 1)\stac L 
c'(k\sut 2 \sw 1 \cdot b')\sigma (k\sut 2 \sw 2, h\sut 2)\stac L
h\sut 1 \sw 3 k\sut 1 \sw 2,
\end{eqnarray*}
where the notation, $\theta^{-1}(h\sstac L 1_\hH)=h\sut 1\sstac {L^{op}} h\sut
2$ (with implicit summation) is used for the inverse of the Galois map
$\theta(h\sstac {L^{op}} k)=\Delta_\hH (h)k$, for $h,k\in \hH$.
\item If $\sigma$ is an invertible cocycle (with inverse $\widetilde \sigma$),
  then $\cD $ is a left $B$-bialgebroid with coproduct and counit
\begin{eqnarray*}
&&\Delta_\cD(b\stac L b'\stac L h)=
\big(b \stac L {\widetilde \sigma}(h\sw 1\sut 2,h\sw 2)\stac L h\sw 1 \sut
1\big) \stac B 
\big(1_B \stac L b'\stac L h\sw 3\big)\qquad \textrm{and}\\
&&\varepsilon_\cD\ (b\stac L b'\stac L h)= 
b(h\sw 1\cdot b')\sigma(h\sw 2 \sut 1 , h\sw 2 \sut 2).
\end{eqnarray*}
\item If $\sigma$ is an invertible cocycle (with inverse $\widetilde \sigma$),
  then $\cD $ is a left $\times_B$-Hopf algebra. That is, the Galois map $\cD
  \sstac {B^{op}} \cD \to \cD \sstac B \cD$, $(b\sstac L b'\sstac L h)\sstac
  {B^{op}} (c\sstac L c'\sstac L k)\mapsto \big(\Delta_\cD (b\sstac L b'\sstac
  L h)\big) (c\sstac L c'\sstac L k)$ is bijective, with inverse 
\begin{eqnarray*}
(b\stac L b'\stac L h)\stac B (c\stac L c'\stac L k)&\mapsto&
\big(b\stac L 1_B \stac L h\sut 1\big) \stac {B^{op}} 
\big(b'(h\sut 2\sw 1\cdot c)\sigma(h\sut 2 \sw 2,k\sut 1\sw 1)\\
&&\stac L 
c' {\widetilde \sigma}(k\sut 2 h\sut 2\sw 3 \sut 2,h\sut 2 \sw 4)\stac L
h\sut 2 \sw 3 \sut 1 k\sut 1\sw 2\big).
\end{eqnarray*}
\end{zlist}
\end{proposition}
Verification of Proposition \ref{prop:twbgd} by a direct (and somewhat
long) computation is left to the reader.
\begin{remark}
(1) Proposition~\ref{prop:twbgd} can be specialised to the case when $\hH$ is
    a Hopf algebra, and $B$ is its module algebra (i.e.\ $\sigma$ is a trivial
    cocycle). Then the $\times_B$-Hopf algebra $\cD$ in the proposition
    reduces to the one constructed in \cite{Kad:Pseudo}, which, subsequently
    has 
    been shown in \cite{PaVOys:K-CMiso} to be isomorphic to the 
    bialgebroid constructed in \cite[Section~3]{ConMos:Ran}.

(2) Apply Proposition \ref{prop:twbgd} to the particular case when the
    $\sigma$-twisted $\hH$-module $B$ is equal to the base algebra $L$
    itself. Then $\cD$ is isomorphic to $\hH$, as a $k$-module. If $\sigma$ is
    an invertible cocycle, then the $k$-linear automorphism of $\hH$,
$
h\mapsto t\big(\sigma(h\sw 2 \sut 1, h\sw 2 \sut 2)\big) h\sw 1,
$ with inverse 
$
h\mapsto h\sw 1 \sut 1 t\big({\widetilde \sigma}(h\sw 1 \sut 2, h\sw 2)\big)
$,
maps the left $\times_B$-Hopf algebra $\cD$ to the {\em cocycle double twist}
of $\hH$. By a cocycle double twist of a left bialgebroid we mean the following
generalisation of Doi's construction on bialgebras in \cite{Doi:BraQu}.
Let $\hH$ be a left $L$-bialgebroid, with source map $s$ and target map $t$.
Assume that $\hH$ measures $L$ and let $\sigma$ be an $L$-valued invertible
2-cocycle. Then $\hH$ is a left $L$-bialgebroid with  
unmodified source and target maps and $L$-coring structure, and a  newly
defined product 
$$
h\stac k h'\mapsto  s\big(\sigma(h\sw 1, h'\sw 1)\big)
t\big({\widetilde \sigma}(h\sw 3, h'\sw 3)\big) h\sw 2 h'\sw 2, \qquad
\textrm{for\ }h,h'\in \hH.
$$
\end{remark}
\begin{theorem}
Let $\hH$ be a Hopf algebroid, with a constituent left $L$-bialgebroid,
right $A$-bialgebroid, and antipode $S$. Let $T$ be a right faithfully flat
$A$-ring and a left 
faithfully flat right $\hH$-cleft extension of $B$. Then the canonical left
$\times_B$-Hopf algebra in Theorem~\ref{thm.main.Hobst}~(iii) 
is isomorphic to $\cD$ in Proposition \ref{prop:twbgd}. 
\end{theorem}
\begin{proof}
Let $j:\hH\to T$ be a normalised cleaving map, with convolution inverse
${\widetilde j}$. Recall from \cite[Theorem 3.12]{BohBrz:cleft} that $T$ is
isomorphic to $B\sstac L \hH$, as a left $B$-module and right
$\hH$-comodule. Hence a $k$-linear isomorphism $(T\sstac A T)^{co\hH}\to \cD$
is given by the map 
\begin{equation}\label{eq:cliso}
u\stac A v \mapsto u\sco 0 \su 0 {\widetilde j}(u\sco 0\su 1)\stac L v\su 0
{\widetilde j}(v\su 1)\stac L u\sco 1,
\end{equation}
with inverse 
$
b\sstac L b'\sstac L h \mapsto b j(h\suc 1)\sstac A b' j\big(S(h\suc 2)\big).
$
Recall that the constituent left $L$-bialgebroid in a Hopf algebroid is a
$\times_L$-Hopf algebra, with inverse Galois map $\theta^{-1}(h\sstac L
k)=h^{(1)} \sstac {L^{op}} S(h^{(2)})k$, for $h\sstac L k\in \hH\sstac L \hH$.
Using the explicit form $u \mapsto u\sco 0\su 0\sstac A {\widetilde j}(u\sco
0\su 1)\sstac B j(u\sco 1)$ of the torsor map, for $u\in T$, and relations
\cite[(4.18) and (4.19)]{BohBrz:cleft} between cleaving maps and 2-cocycles,
the reader can easily check that the isomorphism \eqref{eq:cliso} preserves the
$B$-bialgebroid structure indeed.
\end{proof}

\section{Pre-torsors as differentiable bimodules}
In this appendix we calculate explicitly differential graded algebras arising 
from pre-torsors and describe differentiable bimodule structures on a
pre-torsor.  
Such bimodules are one of the main ingredients  in the construction of
 non-commutative differential fibrations \cite{BegBrz:ser}.
Recall that, given a differential graded algebra 
$(\Omega(A), d)$ with $\Omega^0 (A) =A$ and a right $A$-module $M$, 
a map $\nabla: M\to M\sstac{A}\Omega^1(A)$ is called a (right) {\em connection}
in $M$, provided that for all $a\in A$ and $m\in M$,
$$
\nabla (ma) = \nabla(m)a +m\stac A da.
$$
Similarly, one defines a connection in a left $A$-module. The map $\nabla$ can
be uniquely extended to the map $\nabla: M\sstac{A} \Omega^\bullet(A)\to  
M\sstac{A}\Omega^{\bullet +1}(A)$ by the (graded in the left module case)
Leibniz rule:  
$\nabla(m\sstac{A}\omega)= \nabla(m)\omega + m\sstac{A}d\omega$. 
A connection $\nabla$ is said to be {\em flat}, provided $\nabla\circ\nabla
=0$. 

We say that an $A$-$B$ pre-torsor $T$ is {\em unital} if the structure map
$\tau$ satisfies condition (d) in Definition~\ref{def:ABtor}. For example, 
the pre-torsor constructed in Example~\ref{ex.homog} is unital. 
If, in addition, $T$ is a faithfully
flat pre-torsor, then $A$-coring $\cC$ in
Theorem~\ref{thm.main.1}
has a group-like element $1_T\sstac{B} 1_T$. Now, 
\cite[29.7, 29.14]{BrzWis:cor} yield the following corollary of
Theorem~\ref{thm.main.1}. 
\begin{corollary}\label{prop.con}
Let $T$ be a faithfully flat unital $A$-$B$ pre-torsor with the torsor map
$\tau$. Then 
\begin{zlist}
\item 
\begin{blist}
\item There exists a differential graded (tensor) algebra structure
  $\Omega(A)$ on $A$ with 
$$
\Omega^1(A) := \{\sum_i t_i\stac{B} u_i\in T\stac{B}T \; |\; \sum_i t_iu_i =0\;
\mathrm{and}\;  
\sum _i t_i\tau(u_i) = \sum_i 1_T\stac{A}
t_i\stac{B} u_i \}, 
$$
$
\Omega^n(A) = \Omega^1(A)^{\sstac{A}n},
$
and differentials
$
d(a) = 1_T\sstac{B}\alpha(a) - \alpha(a)\sstac{B} 1_T$, for all $a\in A$ and,
for all $\sum_i t_i\sstac{B} u_i\in \Omega^1(A)$, 
$$
d(\sum_i t_i\stac{B} u_i) = \sum_i 1_T\stac{B} 1_T\stac{A} t_i\stac{B} u_i
- \sum_i t_i\stac{B} \tau(u_i) + \sum_i t_i\stac{B} u_i\stac{A}1_T\stac{B} 1_T.
$$
\item The map 
$$
\nabla_T^r: T\to T\stac{A}\Omega^1(A), \qquad 
t\mapsto \tau(t) - t\stac{A}1_T\stac{B}1_T,
$$
is a flat connection in $T$.
\end{blist}
\item 
\begin{blist}
\item There exists a differential graded (tensor) algebra structure
  $\Omega(B)$  on $B$ with 
$$
\Omega^1(B) := \{\sum_i t_i\stac{A} u_i \in T\stac{A}T\; |\; \sum_i t_iu_i =0\;
\mathrm{and}\;  
\sum _i \tau(t_i)u_i =\sum_i t_i\stac{A} u_i \stac{B} 1_T\}, 
$$
$
\Omega^n(B) = \Omega^1(B)^{\sstac{B}n},
$
and differentials
$
d(b) = 1_T\sstac{A}\beta(b) - \beta(b)\sstac{A} 1_T$, for all $b\in B$ and,
for all $\sum_i t_i\sstac{A} u_i\in \Omega^1(B)$, 
$$
d(\sum_i t_i\stac{A} u_i) = \sum_i 1_T\stac{A} 1_T\stac{B} t_i\stac{A} u_i
- \sum_i \tau(t_i)\stac{A} u_i + \sum_i t_i\stac{A} u_i\stac{B}1_T\stac{A} 1_T.
$$
\item The map 
$$
\nabla_T^l: T\to \Omega^1(B)\stac{B}T, \qquad 
t\mapsto  1_T\stac{A}1_T\stac{B} t -\tau(t), 
$$
is a flat connection in $T$.
\end{blist}
\end{zlist}
\end{corollary}

Recall from \cite[Definition~2.10]{BegBrz:ser} (cf.\
 \cite[Section~3.6]{Mad:int}) 
 that, given 
differential graded algebras $\Omega(A)$ over $A$ and $\Omega(B)$ over $B$, a 
$B$-$A$-bimodule $M$ with a left connection 
$\nabla: M\to \Omega^1(B)\sstac{B} M$ and a $B$-$A$ bimodule map 
$\sigma: M\sstac{A}\Omega^1(A)\to \Omega^1(B)\sstac{B}M$ is called a 
{\em (left) differentiable bimodule}, provided that, for all $m\in M$ and
$a\in A$, $\nabla(ma) = \nabla(m)a +\sigma(m\sstac{A}da)$. In case $A=B$,
the pair $(\nabla,\sigma)$ is termed a {\em (left) $A$-bimodule connection} 
\cite{DubMad:cur}. 
Differentiable
bimodules induce functors between categories of connections (cf.\ 
\cite[Proposition~2.12]{BegBrz:ser}). Since the left connection
$\nabla_T^l$ in 
Corollary~\ref{prop.con}(2)(b) is right $A$-linear, every
faithfully flat unital $A$-$B$ torsor is a differentiable bimodule with
the trivial (zero) map $\sigma$. The corresponding functor coincides with
the tensor (induction) functor $T\sstac{A}\, \bullet$. More interesting
 differential 
bimodule structure can be constructed on $T$, provided the structure map
$\tau$ is $B$-$B$ bilinear.

\begin{proposition}\label{prop.dif.mod}
Let $T$ be a faithfully flat unital $A$-$B$ pre-torsor with differential
structures $\Omega(A)$ and  
 $\Omega(B)$ as in Corollary~\ref{prop.con}. 
\begin{zlist}
\item The pair $(\nabla_T^l,\sigma_B)$, where
$$
\sigma_B: T\stac{B}\Omega^1(B)\to \Omega^1(B)\stac{B}T, \qquad 
\sum_i t\stac{B} u_i\stac{A} v_i \mapsto \sum_i \tau(t u_i)v_i.
$$
 is a $B$-bimodule connection.
 \item If the torsor
structure map $\tau$ is right $B$-linear, then $T$ is a 
differentiable $B$-$A$ bimodule with (flat) connection $\nabla_T^l$ in 
Corollary~\ref{prop.con}(2)(b) and the twist map
$$
\sigma^l: T\stac{A}\Omega^1(A)\to \Omega^1(B)\stac{B}T, \qquad 
\sum_i t\stac{A} u_i\stac{B} v_i \mapsto \sum_i \tau(t u_i)v_i.
$$
\end{zlist}
\end{proposition}
\begin{proof}
(1) 
First we need to check whether the map $\sigma_B$ is well-defined. The
property (a) in Definition~\ref{def.pre.tor} implies that
$\im \sigma_B\subseteq \ker(\mu_T\sstac{B}T)$.
Second, by properties
(c) and (a) in Definition~\ref{def.pre.tor}, 
$$(T\stac{A}T\stac{B}\mu_T\stac{B} T)\circ (\tau\stac{A}T\stac{B}T)\circ
\sigma_B = (T\stac{A}T\stac{B}\beta\stac{B} T)\circ\sigma_B.
$$
Since $T$ is a (faithfully) flat left $B$-module, this implies that the map
$\sigma_B$ 
is well-defined. An easy calculation verifies the twisted Leibniz rule.

 (2) The map $\sigma^l$ is 
well-defined by the same arguments as in part (1) and the right $B$-linearity
of $\tau$. 
Since $\nabla^l_T$ is a right $A$-module map,  we 
only need to show that, for all $a\in A$ and $t\in T$, $\sigma^l(t\sstac{A}da)
= 0$. This 
follows by the right $A$-linearity of $\tau$. 
\end{proof}

Note that the twisting map $\sigma_B$ in Proposition~\ref{prop.dif.mod} is a
restriction of $\psi_\cD$ in \eqref{eq:psiD}.

\section*{Acknowledgements} 
The work of the first author is supported by the Hungarian Scientific 
Research Fund OTKA T 043 159 and the Bolyai J\'anos Fellowship. She thanks
C. Menini for useful discussions.
The second author would like to express his gratitude to the members of the 
Department of Theoretical Physics at the Research Institute for Particle and
Nuclear Physics in Budapest for very warm hospitality.


\begin{thebibliography}{Bibliography}{} 
\bibitem{Abu:Mor.cor} J. Abuhlail, {\em Morita contexts for corings and
  equivalences}, in S. Caenepeel and F. Oystaeyen (eds.)
  Lecture Notes in Pure and Applied Math.\ Vol. 239 Marcel Dekker New York
  2005, pp 1--19.
\bibitem{BegBrz:ser} E.J.\ Beggs and T.\ Brzezi\'nski, {\em The Serre spectral 
sequence of a noncommutative fibration for de Rham cohomology}, Acta Math.\ 195
(2005), 155--196.
\bibitem{Bohm:hgdint} G.\ B\"ohm, {\em Integral theory for Hopf algebroids},
Alg.\ Rep.\ Theory 8(4) (2005), 563--599. 
\bibitem{BohBrz:str}  G.\ B\"ohm and T.\ Brzezi\'nski, {\em Strong connections
  and
the relative Chern-Galois character for corings},
 Int.\ Math.\ Res.\ Notices, 2005:42 (2005), 2579--2625.
 \bibitem{BohBrz:cleft} G.\ B\"ohm and T.\ Brzezi\'nski, {\em Cleft extensions
  of Hopf algebroids}, 
  Appl.\ Cat.\ Str.\ 14  (2006), 431-469.
\bibitem{BohmSzl:hgdax} G.\ B\"ohm and K.\ Szlach\'anyi, {\em Hopf
algebroids with bijective antipodes: axioms, integrals and duals},
J.\ Algebra {274} (2004), 585-617. 
\bibitem{BohmVer:Mor&cleft} G. B\"ohm and J. Vercruysse, {\em Morita theory
  for coring extensions and cleft bicomodules}, Adv.\ Math.\ 209 (2007),
  611--648.
\bibitem{Brz:corext} T.\ Brzezi\'nski, {\it A note on coring
extensions,} Annali dell'Universit\'a di Ferrara, Sezione VII.-Sci.\ Mat.\ 51
  (2005), 15--27.
\bibitem{BrzHaj:coa}
T.~Brzezi\'nski and  P.M.~Hajac, {\em Coalgebra extensions and 
algebra coextensions of
Galois type},  Comm.\ Alg.\ 27 (1999), 1347--1367.
\bibitem{BrzWis:cor} T.\ Brzezi\'nski and R.\ Wisbauer, {\em Corings and 
 Comodules}. Cambridge University Press, Cambridge, 2003. 
 Erratum: http://www-maths.swan.ac.uk/staff/tb/corinerr.pdf
\bibitem{ConMos:Ran} A.\ Connes and H.\ Moscovici, {\em 
Rankin-Cohen brackets and the Hopf algebra of transverse geometry},  
Mosc.\ Math.\ J.\ 4 (2004), 111--130.
\bibitem{Doi:BraQu} Y.\ Doi, {\em Braided bialgebras and quadratic bialgebras},
  Comm. Algebra 21 (1993), 1731--1749.
\bibitem{DelMil:TanCat} P. Deligne and J. Milne, {\em Tannakian Categories},
  Lecture Notes Math. 900, Springer 1982, pp 101--228. 
\bibitem{DubMad:cur} M.\ Dubois-Violette, J.\ Madore, T.\ Masson, J.\ Mourad,
{\em On curvature in noncommutative geometry}, J.\ Math.\ Phys.\ 37 (1996), 
4089--4102.
\bibitem{ElKaGomTorLob:Semisimp} L. El Kaoutit, J. G\'omez-Torrecillas and
  F.J. Lobillo, {\em Semisimple corings}, Algebra Colloquium 11:4 (2004),
  427--442. 
\bibitem{Gru:tor} C.\ Grunspan, {\em Quantum torsors}, J.\ Pure Appl.\ Alg.\ 
184 (2003), 229--255.
\bibitem{Hobst:PhD} D.\ Hobst, {\em Antipodes in the theory of noncommutative
  torsors}, PhD thesis Ludwig-Maximilians Universit\"at M\"unchen (2004),
  Logos Verlag Berlin, 2004.  
\bibitem{Kad:Pseudo} L.~Kadison, {\em Pseudo-Galois extensions and Hopf
  algebroids}, arXiv:math.QA/0508411.
\bibitem{KadSzl:D2bgd} L.\ Kadison and K.\ Szlach\'anyi, {\em
Bialgebroid actions on depth two extensions and duality}, Adv.\
Math.\ {179} (2003), 75-121.    
\bibitem{Kon:ope} M.\ Kontsevich, {\em Operads and motives in deformation
  quantisation}, Lett.\ Math.\ Phys.\ 48 (1999), 35--72. 
\bibitem{Lu:hgd} J.-H.\ Lu, {\em Hopf algebroids and quantum groupoids},
Int. J. Math. {7} (1996) 47-70.
\bibitem{Mad:int}
J.\ Madore,
{\em An Introduction to Noncommutative Differential Geometry
    and Its Physical Applications.} 2nd ed., Cambridge
University Press, Cambridge (1999).
\bibitem{MoeMrc:fol} I.\ Moerdijk and J.\ Mr\v cun, {\em Introduction to
  Foliations and Lie Groupoids}, Cambridge University Press, Cambridge (2003).
  \bibitem{PaVOys:K-CMiso} F. Panaite and F. Van Oystaeyen, {\em Some
  bialgebroids constructed by Kadison and Connes-Moscovici are isomorphic},
  Appl.\ Cat.\ Str.\ 14 (2006), 627-632.
\bibitem{MulSch:QHomSpff} E.F.\ M\"uller and H.-J.\ Schneider, Quantum
  homogenous spaces with faithfully flat module structures, {\sl Israel
  J. Math.} {\bf 111} (1999), 157--190.
\bibitem{SaaR:CatTan} N. Saavedra Rivano, {\em Cat\'egories Tannakiennes}
  Lecture Notes in Math. 265, Springer 1972.
\bibitem{Sch:big} P.\ Schauenburg, {\em Hopf bi-Galois extensions},
Commun.\ Algebra 24 (1996), 3797--3825.
\bibitem{Sch:Bianc}  P.\ Schauenburg, {\em Bialgebras over noncommutative
  rings and a structure theorem for Hopf bimodules}, Appl.\ Cat.\
  Str.\ 6 (1998), 193--222.
\bibitem{Sch:dua} P.\ Schauenburg, {\em Duals and doubles of quantum groupoids
 ($\times_R$-Hopf algebras)},
  in N.\ Andruskiewitsch, W.R.\ Ferrer Santos, H.-J.\ Schneider (eds.)
{\em New trends in Hopf algebra theory}, Contemp.\ Math.\ vol.\ 267, 
Amer.\ Math.\ Soc., Providence, RI, (2000), pp.\ 273--299.
\bibitem{Sch:tor} P.\ Schauenburg, {\em Quantum torsors with fewer axioms},
arXiv: math.QA/0302003. 
\bibitem{Sch:Hop} P.\ Schauenburg, {\em Hopf-Galois and bi-Galois extensions},
  in 
{\em Galois theory, Hopf algebras, and semiabelian categories}, Fields Inst.\
  Commun., 43, Amer.\ Math.\ Soc., Providence, RI, (2004), pp.\ 469--515. 
\bibitem{Schn:PriHomSp} H.-J. Schneider, Principal homogenous spaces for
  arbitrary Hopf algebras, {\sl Israel J. Math.} {\bf 72} (1990), 167--195.
\bibitem{Tak:crR} M.\ Takeuchi, {\em Groups of algebras over $A\ot {\bar
A}$}, J. Math. Soc. Japan {29} (1977), 459-492.
\bibitem{Ulb:Galfun} K.-H. Ulbrich, {\em Galois extensions as functors of
  comodules}, Manuscripta Math. 59 (1987), 391--397.
\bibitem{Ulb:Fibre} K.-H. Ulbrich, {\em Fibre functors of finite dimensional
  comodules}, Manuscripta Math. 65 (1989), 39--46.
\bibitem{Ulb:Rec} K.-H. Ulbrich, {\em On Hopf algebras and rigid monoidal
  categories}, Israel J. Math. 72 (1990), 252--256.
\end{thebibliography}
\end{document}